\documentclass[12pt]{article}

\usepackage{amsmath, amssymb, amsthm}
\usepackage{graphicx}

\newcommand{\e}{\textbf{e}}
\newcommand{\et}{\tilde{\textbf{e}}}

\newcommand{\PT}{\mathcal{PT}}
\newcommand{\as}{$\alpha$-surface}
\newcommand{\bs}{$\beta$-surface}
\newcommand{\OO}{\mathcal{O}}
\newcommand{\CP}{\mathbb{CP}^1}
\newcommand{\RP}{\mathbb{RP}^1}
\newcommand{\cK}{\mathcal{K}}
\newcommand{\tK}{\tilde{K}}

\newcommand{\ep}{\epsilon}
\newcommand{\cL}{\mathcal{L}}

\newcommand{\D}{\mathcal{D}}

\newcommand{\R}{\mathbb{R}}
\newcommand{\C}{\mathbb{C}}
\newcommand{\cS}{\Theta}

\renewcommand{\d}{\mathrm{d}}
\newcommand{\bK}{\mathbb{K}}
\newcommand{\p}{\partial}

\newcommand{\lam}{\lambda}

\newcommand{\ot}{\otimes}

\newcommand{\Sm}{\mathbf{\Sigma}}
\long\def\symbolfootnote[#1]#2{\begingroup%
\def\thefootnote{\fnsymbol{footnote}}\footnote[#1]{#2}\endgroup}

\newtheorem{theorem}{Theorem}
\newtheorem{proposition}{Proposition}
\newtheorem{lemma}{Lemma}
\newtheorem*{thm}{Theorem}

\def\be{\begin{equation}}
\def\ee{\end{equation}}

\title{Anti-self-dual conformal structures with null Killing vectors
  from projective structures}
\author{Maciej Dunajski, Simon West}

\begin{document}

\pagestyle{plain}

\title{\vskip -70pt
\begin{flushright}
{\normalsize DAMTP-2006-9} \\
\end{flushright}
\vskip 50pt
{\bf Anti-self-dual conformal structures with null Killing vectors
  from projective structures} \vskip 20pt}

\author{Maciej Dunajski\thanks{email M.Dunajski@damtp.cam.ac.uk} \ \
and \  \
Simon West\thanks{email S.D.C.West@damtp.cam.ac.uk} \\[15pt]
{\sl Department of Applied Mathematics and Theoretical Physics} \\[5pt]
{\sl University of Cambridge} \\[5pt]
{\sl Wilberforce Road, Cambridge CB3 0WA, UK} \\[15pt]}
\date{January 17, 2006} 
\maketitle
\begin{center} \emph{Dedicated to the memory of Jerzy Pleba\'nski} \end{center}
\begin{abstract}
Using twistor methods, we explicitly construct all local forms of 
four--dimensional real analytic neutral signature 
anti--self--dual conformal structures
$(M,[g])$ with a null
conformal Killing vector. We show that $M$ is foliated by
anti-self-dual null surfaces, and the two-dimensional leaf space 
inherits a natural projective structure. The twistor space
of this projective structure is the quotient of 
the twistor space of $(M,[g])$ by the group action induced by the 
conformal Killing vector.

We obtain a local classification  which
branches according to whether or not the conformal Killing  vector 
is hyper-surface orthogonal in $(M, [g])$.  We give examples of conformal 
classes which contain Ricci--flat metrics on compact complex surfaces
and discuss other conformal classes with no Ricci--flat metrics.
\end{abstract}
\newpage
\section{Introduction}
\setcounter{equation}{0}
The anti--self--duality (ASD) condition in four dimensions
seems to underlie the concept of integrability of ordinary and partial
differential equations \cite{Wa85}. Many lower dimensional integrable models
(KdV, NlS, Sine--Gordon, ...)
arise as symmetry reductions of the ASD Yang--Mills equations on a flat
background, and various solution generation techniques are reductions
of the twistor correspondence \cite{masonwoodhouse}.
Other integrable models (dispersionless Kadomtsev--Petviashvili, $SU(\infty)$ Toda, ...) 
are reductions of the ASD conformal equations
which say that the self--dual Weyl tensor of a conformal class of metrics
vanishes \cite{W90,DMT}. 
Generalisations to ASD Yang--Mills on ASD conformal background are also 
possible \cite{Tafel0, calderbank}.

In all cases the main interest is in conformal structures
of signature $++--$ which are called \emph{neutral}, as the
reductions can lead to interesting  hyperbolic and parabolic equations. 
There are no non--trivial ASD structures in the Lorentzian signature
$+++-$, and the reductions from Riemannian manifolds can
only yield elliptic equations thus ruling out
interesting soliton dynamics.

The main gap in the programme to classify the reductions of 
ASD neutral conformal structures was  understanding the reductions
by a null conformal Killing vector. 
We embarked on this project
hoping to incorporate more integrable systems into the framework of 
of  anti--self--duality, but we have found (Theorem 2) that the
resulting geometry is a  \emph{completely solvable} system.

Let $(M,[g])$ be a four dimensional real analytic
neutral ASD conformal structure. We say that $K$ is a 
null conformal Killing vector if it satisfies 
\begin{equation}
\label{definition_K}
\mathfrak{L}_K g = \eta g, \qquad g(K, K)=0,
\end{equation}
for some $g\in [g]$, where $\eta$ is a function on $M$, and 
$\mathfrak{L}$ is the Lie-derivative.

When studying conformal structures with non-null conformal Killing
vectors, it is natural to look at the space of Killing vector
trajectories, since this will inherit a non-degenerate conformal
structure. In the case of a \emph{null} conformal Killing vector, the
situation is different. The space of trajectories inherits a
\emph{degenerate} conformal structure. We find that it is necessary to 
go down one
dimension more, and consider a \emph{two} dimensional space $U$ of
anti-self-dual totally null surfaces in $M$, called $\beta$--surfaces, 
containing
$K$, which exist as a consequence of the conformal Killing
equation. It turns out that there is a naturally defined projective
structure $[\Gamma]$ on $U$. Moreover, we show that the twistor
spaces of $(M,[g])$ and $(U,[\Gamma])$ are related by dimensional
reduction. Specifically, the twistor space $Z$ of $(U,[\Gamma])$ is
the space of trajectories of a vector field on the twistor space
$\PT$ of $(M,[g])$ corresponding to $K$.
Projective structures are just equivalence classes of torsion-free connections,
which do not need to satisfy any equations; 
this underlies the complete solvability of null reductions,
and contrasts with the
non-null case where one obtains Einstein-Weyl structures
\cite{jonestod}, and associated integrable systems
\cite{W90,DMT,dunajski,calderbank}.

In Section \ref{NullKVs} we
derive some elementary properties of null conformal Killing
vectors. Section \ref{projstructures} is an introduction to projective
structures. In Section \ref{main} we prove the following:
\begin{theorem} \label{maintheorem}
Let $(M,[g])$ be a four dimensional real analytic
neutral ASD conformal structure with a null conformal Killing
vector $K$. Let $U$
be the two dimensional space of \bs s containing $K$. Then there is a
naturally defined projective structure on $U$, whose twistor space is
the space of trajectories of a distribution $\widehat{\cK}$ induced on
$\PT$ by the action of $K$ on $M$.
\end{theorem}

In Section \ref{local} we investigate the local form of ASD conformal
structures with null Killing vectors. This is expressed in the
following theorem:

\begin{theorem} \label{localform}
Let $(M,[g],K)$ be a smooth neutral signature ASD conformal structure
with null conformal Killing vector. Then there exist
local coordinates $(t,x,y,z)$ and $g \in [g]$ such that $K=\p_t$ and  
$g$ has one
of the following two forms, according to whether the  twist
$\bK \wedge d \bK$
vanishes or not $(\bK := g(K,.))$:
\begin{enumerate}
\item $\bK \wedge d \bK =0$.
\begin{multline} \label{nontwistinggeneral}
g = (dt + (z A_3 - Q) dy)(dy - \beta dx) - \\(dz - (z (-\beta_y +
A_1 + \beta A_2 + \beta^2 A_3))dx - (z(A_2 + 2 \beta A_3)+ P)dy)dx,
\end{multline}
where $A_1,A_2,A_3,\beta,Q,P$ are arbitrary functions of $(x,y)$.

\item $\bK \wedge d \bK \neq 0$.
\begin{multline} \label{twistinggeneral}
g = (dt + A_3\p_zG  dy + (A_2\p_z G  + 2 A_3 (z \p_zG - G)
- \p_z\p_yG)dx)(dy - z dx) \\  - \p_z^2G dx(dz - (A_0 + z A_1 + z^2
A_2 + z^3 A_3)dx),
\end{multline}
where $A_0,A_1,A_2,A_3$ are arbitrary functions of $(x,y)$, and
$G$ is a function of $(x,y,z)$ satisfying the following PDE:
\begin{equation} \label{Gequation}
(\p_x + z \p_y + (A_0+zA_1+z^2 A_2 + z^3 A_3) \p_z) \p_z ^2G = 0.
\end{equation}
\end{enumerate}
\end{theorem}
The functions $A_\alpha(x,y)$ in the metrics (\ref{nontwistinggeneral}) and
(\ref{twistinggeneral}) determine projective structures on the
two dimensional space $U$ in the following way. A general projective structure
corresponds to a second-order ODE
\begin{equation} \label{secondorderode}
\frac{d^2 y}{dx^2} = A_3(x,y) \Big( \frac{d y}{dx} \Big)^3 +
A_2(x,y) \Big( \frac{dy}{dx} \Big)^2 + A_1(x,y) \Big( \frac{dy}{dx}
\Big) + A_0(x,y).
\end{equation}
In (\ref{twistinggeneral}) all the $A_\alpha, \alpha=0,1,2,3$ functions occur
explicitly in the metric. In (\ref{nontwistinggeneral}) the function
$A_0$ does not explicitly occur. It is determined by the following
equation:
\be
\label{beta_equation}
A_0 = \beta_x + \beta \beta_y - \beta A_1 - \beta^2 A_2 - \beta^3 A_3,
\ee
as is shown in the proof of the theorem. 
If one interprets
$z$ as a fibre coordinate on the
projective tangent bundle of the $(x,y)$ space, then
(\ref{Gequation}) says
that $\p_z^2 G$ is constant along the projective structure
spray, (compare formula \ref{nonhomspray}). 

Note that in both cases the Killing vector is $\p_t$ and is pure
Killing (this comes from choosing a suitable $g\in [g]$).
The non-twisting case (\ref{nontwistinggeneral})
is a
natural conformal generalisation of Ricci--flat pp waves.
The twisting case (\ref{twistinggeneral}) is a neutral analog
of the  Fefferman conformal class \cite{graham}. As special cases of
(\ref{twistinggeneral}) we
recover some examples of \cite{nurowskisparling}, where neutral
metrics were related to second order ODEs.

The aim of Section \ref{examples} is to put our results into a broader
context.
We examine some examples found by different
means in the light of our results. %Section \ref{specialmetrics} is
%concerned with examining the local conformal structures in more
%detail. For example we find Ricci-flat metrics, and also examples of
%conformal classes containing no Ricci-flat metrics.
We find necessary and sufficient conditions on the underlying
projective structure in order for there to exist (pseudo)
hyper--complex metrics with triholomorphic $K$ within a conformal class.
A special case of the metric (\ref{nontwistinggeneral})  yields a
compact example of a Ricci--flat metric on a Kodaira surface of a special type.

We consider how to construct conformal structure twistor spaces from
projective structure twistor spaces in Section \ref{twistorreconstruction}. 
The more involved spinor calculations are moved to the Appendix. 

\section{Null Conformal Killing Vectors} \label{NullKVs}
\setcounter{equation}{0}
\subsection{Spinors in neutral signature} \label{spinorsinneutralsig}

We will denote by $(M,[g])$ a local patch of $\R^4$ endowed with a
neutral signature conformal structure $[g]$. That is, $[g]$ is an
equivalence class of neutral signature metrics with the equivalence
relation $g \sim e^c g$ for a some function $c$ on $M$.

Any neutral metric $g$ on $M$ can be put in the following form:
\begin{equation} \label{tetradform}
g=2(\mathbf{\theta}^{00'} \odot \mathbf{\theta}^{11'} -
\mathbf{\theta}^{01'} \odot \mathbf{\theta}^{10'}) =
\ep_{AB}\ep_{A'B'}\mathbf{\theta}^{AA'} \otimes \mathbf{\theta}^{BB'},
\end{equation}
where $\ep_{AB}, \ep_{A'B'}$ are antisymmetric matrices with $\ep_{01}
= \ep_{0'1'} = 1$. The
four (real) basis one-forms $\mathbf{\theta}^{AA'}$ for $A=0,1$,
$A'=0,1$ are called a tetrad. The algebraic dual vector basis is
denoted $\e_{AA'}$, and is defined by
$\mathbf{\theta}^{AA'}(\e_{BB'})=\delta_B^A \delta_{B'}^{A'}$. Any vector $V$ at a
point can be written $V^{AA'}\e_{AA'}$, and this exhibits an
isomorpism
\begin{equation} \label{spinbundles}
TM \cong S \otimes S',
\end{equation}
where $S$, $S'$ are two-dimensional real vector bundles known as the
\emph{unprimed spin bundle} and the \emph{primed spin bundle}
respectively. For a general manifold $M$ there is a topological
obstruction to ($\ref{spinbundles}$) but we are working locally so it
always holds.

Using a particular choice of tetrad, a section $\mu$ of $S$ is denoted
$\mu^{A}$, $A=0,1$. Similarly $\nu_{A}$ is a section of $S^*$,
$\kappa^{A'}$ a section of $S'$ and $\tau_{A'}$ a section of
$S'^*$, where $^*$ denotes the dual of a bundle. The natural pairing
$S \times S^* \rightarrow \mathbb{R}$ is
given by $\mu^A \nu_A$, using the summation convention, and similarly
for primed spinors. We sometimes use the notation $\mu^{A'} \nu_{A'} =
\mu.\nu$. This product is not commutative, we have $\mu.\nu = -\nu.\mu$.

It follows from (\ref{tetradform}) that $g(V,V)=\text{det} \
V^{AA'}$. If $V$ is null, then this gives $V^{AA'} = \mu^A
\kappa^{A'}$. Abstractly, if $V$ is null then $V=\mu \otimes \kappa$ under
the isomorphism
(\ref{spinbundles}), where $\mu$ and $\kappa$ are sections of $S$ and
$S'$ respectively.

The relation (\ref{tetradform}) can be written abstractly as
$$
g = \ep \otimes \ep'
$$
under the isomorphism $(\ref{spinbundles})$. $\ep$ and $\ep'$ are
symplectic structures on $S$ and $S'$. These give isomorphisms $S
\cong S^*$ and $S'
\cong S'^*$ by $\mu \rightarrow \ep(\mu,.)$, for $\mu$ a section of
$S$, and similarly for $S'$. Given a choice of tetrad, the spinors
$\ep$ and $\ep'$ are written $\ep_{AB}$ and $\ep_{A'B'}$, where we drop the prime
on the latter because no confusion can arise due to the indices. Note
these are anti-symmetric in $AB$ and $A'B'$. Then
the isomorphism $S \cong S^*$ is given in the trivialization by
$\mu^{A} \rightarrow \mu^B \ep_{BA} := \mu_A$ and similarly for primed spinors.

There are useful isomorphisms
\begin{equation} \label{twoformdecomp}
\Lambda_+^2 \cong \text{Sym}(S^* \otimes S^*), \ \ \ \Lambda_-^2 \cong
\text{Sym} (S'^* \otimes S'^*),
\end{equation}
where $\Lambda_+$, $\Lambda_-$ are the bundles of self-dual and
anti-self-dual two-forms, using an appropriate choice of volume form
for the Hodge-$\ast$ operator.
In the local trivialization, the isomorphisms (\ref{twoformdecomp}) are
expressed by the following formula for a two-form $F$ in spinors:
$$
F_{ab} = F_{AA'BB'} = \phi_{A'B'} \epsilon_{AB} + \psi_{AB} \epsilon_{A'B'},
$$
where $\phi_{A'B'}$, $\psi_{AB}$ are symmetric. The $\phi_{A'B'}$ term
is the self-dual component of $F$ and the $\psi_{AB}$ is the
anti-self-dual component.

The vector bundles $S$, $S'$ and their duals inherit connections from
the Levi-Civita connection of $TM$ (see Appendix \ref{appendix}). These are
the unique torsion-free connections defined so that the sections
$\ep$ and ${\ep}'$ are covariantly constant. Then covariant
differentiation on either side of $(\ref{spinbundles})$ is
consistent.

A primed spinor $\kappa^{A'}$ at a point
corresponds to a totally null self-dual two-plane spanned by $\kappa^{A'} \e_{AA'}$,
$A=1,2$, whilst an unprimed spinor corresponds to an anti-self-dual
two-plane in a similar way. In twistor theory, these two-planes are called
$\alpha$-planes and $\beta$-planes respectively. 

\subsection{Null conformal Killing vectors in neutral signature} \label{gsf}

Suppose $g$ is a neutral metric with a conformal Killing
vector $K$. Then $\mathfrak{L}_K (e^c g) = (K(e^c)+e^c \eta)g$, so
$K$ is
a conformal Killing vector for the conformally rescaled metric, and we
can refer to $K$ as a conformal Killing vector for the conformal
structure $[g]$.

Now suppose $g$ has a null conformal Killing vector $K$. We shall show
(Lemma \ref{kvproperties}) that
$M$ is foliated in two different ways, by self-dual and anti-self-dual
surfaces, whose leaves intersect tangent to $K$. This is a property of
the conformal structure $[g]$, since the Hodge-$\ast$ acting on 2-forms
is conformally invariant.

The spinor form of the conformal Killing equation is:
\begin{equation}\label{killing}
\nabla_{a}K_{b}=\phi_{A'B'} \epsilon_{AB} + \psi_{AB} \epsilon_{A'B'}
+ \frac{1}{2} \eta \ep_{AB}\ep_{A'B'},
\end{equation}
where $\phi_{A'B'}, \psi_{AB}$ are the self-dual and anti-self dual
parts of the 2-form $\nabla_{[a}K_{b]}$, and $\eta$ is a function on $M$.

Since $K$ is null, we have $K = \iota \otimes o$, where $\iota$ is a
section of $S$ and $o$ a section of $S'$. Choosing a null tetrad, 
and  a trivialization of $S$ and $S'$,
we have $K^{AA'} = \iota^A o^{A'}$. These spinors are
defined up to multiplication by a non-zero function $\alpha$, since
$K^{AA'}=\iota^A o^{A'}=(\alpha \iota^A)(o^{A'} / \alpha)$.
\begin{lemma} \label{kvproperties}
Let $K=\iota^A o^{A'} \e_{AA'}$ be a null conformal Killing vector. Then 
\begin{enumerate}
\item The following algebraic identities hold:
\begin{eqnarray}
\iota^A \iota^B \psi_{AB} = 0, \label{alg1}\\
o^{A'} o^{B'} \phi_{A'B'} = 0. \label{alg2}
\end{eqnarray}
\item $\iota^A$ and $o^{A'}$ satisfy
\begin{eqnarray} 
\iota^A \iota^B \nabla_{BB'} \iota_A = 0, \label{GSF1}\\
o^{A'} o^{B'} \nabla_{BB'} o_{A'} = 0. \label{GSF2}
\end{eqnarray}
\end{enumerate}
\end{lemma}
\textbf{Remark.} The equations (\ref{GSF1}), (\ref{GSF2}) are
equivalent to the statement that the distributions spanned by $\iota^A
\e_{AA'}$ and $o^{A'} \e_{AA'}$ are Frobenius integrable (see Appendix). Equations of this type are often called
`geodesic shear free' equations, since in the Lorentzian case they
result in shear-free congruences of null geodesics.

\vspace{12pt}

\emph{Proof}. Using $K_{AA'}=\iota_A o_{A'}$, the Killing equation
(\ref{killing}) becomes 
\begin{equation} \label{nullkv}
o_{A'} \nabla_{BB'} \iota_A + \iota_A \nabla_{BB'} o_{A'} =
\phi_{A'B'} \epsilon_{AB} + \psi_{AB} \epsilon_{A'B'} +
\frac{1}{2}\eta \ep_{AB}\ep_{A'B'}.
\end{equation}
Contracting both sides with $\iota^A o^{A'}$ gives 
$$
0= o^{A'} \iota_{B} \phi_{A'B'} + \iota^A o_{B'} \psi_{AB} +
\frac{1}{2} \eta \iota_B o_{B'}.
$$
Multiplying by $\iota^B$ and $o^{B'}$ respectively
leads to (\ref{alg1}) and (\ref{alg2}). To get
(\ref{GSF1}) and (\ref{GSF2}), multiply (\ref{nullkv}) by $\iota^A
\iota^B$ and $o^{A'} o^{B'}$. $\square$
\vspace{12pt}

%Note that Lemma \ref{kvproperties} holds for a holomorphic conformal
%structure with null holomorphic conformal Killing vector The reason
%all the results carry through is that the spinor arguments
%apply to both situations. The only difference is that in the real
%split signature case spinors can be real, whereas in the holomorphic case
%they are complex.

We have found that $M$ is foliated in two different ways by totally
null surfaces. Those determined by $o^{A'}$ are self-dual and are
called \as s, and those determined by $\iota^A$ are anti-self-dual and
are called \bs s. It is
clear that the \as s and \bs s of Lemma \ref{kvproperties} intersect
on integral curves of $K$. Denote the \bs \
distribution by $\D_{\beta}$; this will be used later.

It is appropriate here to recall the Petrov-Penrose classification
\cite{penroserindler} of the algebraic type of a Weyl tensor. In split
signature this applies separately to $C_{ABCD}$ and
$C_{A'B'C'D'}$. In our case $C_{A'B'C'D'}=0$ and we are concerned with the
algebraic type of $C_{ABCD}$. When we refer to the algebraic type we
will be referring to the algebraic type of $C_{ABCD}$. One can form a
real polynomial of fourth
order $P(x)$ by defining $\mu^A = (1,x)$ and setting $P(x) = \mu^A
\mu^B \mu^C \mu^D C_{ABCD}$. The Petrov-Penrose classification refers
to the position of roots of this polynomial, for example if there are
four repeated roots then we say $C_{ABCD}$ is type N. If there is a
repeated root the metric is called \emph{algebraically special}.
There are additional complications in the split signature case
\cite{Law} arising from the fact that real polynomials may not have
real roots.

The split signature version of the Goldberg-Sachs theorem together
with (\ref{GSF1}) implies that any Ricci-flat ASD space with null conformal
Killing vector is algebraically special. In fact the vacuum condition
can be removed if $K$ is non-twisting; we will discuss this further in
Section \ref{confvac}.

It also follows from the Killing equations and the fact that $K$ is null that
$$
K^b \nabla_b K_a = \frac{1}{2} \eta K_a.
$$
Thus $K$ is automatically geodesic, and if it is pure then its
trajectories are parameterized by an affine parameter.

\section{Projective structures} \label{projstructures}
\setcounter{equation}{0}
Let $(U,[\Gamma])$ be a local two dimensional real projective
structure. That is, $U$ is a local patch of $\mathbb{R}^2$, and
$[\Gamma]$ is an equivalence class of torsion-free connections whose
unparameterized geodesics are the same. Then in a local
trivialization, equivalent
torsion-free connections are related in the following way:
\begin{equation} \label{connectiondifference}
\tilde{\Gamma}^i_{jk} - \Gamma^i_{jk} = a_j \delta^i_k + a_k \delta^i_j,
\end{equation}
for functions $a_i$ on $U$, and $i,j,k=1,2$. Note that this is a
tensor equation since
the difference between two connections is a tensor. The $a_i$ on the
RHS are the components of a one-form.

The geodesics satisfy the following ODE:
$$
\frac{d^2 s^i}{dt^2}+\Gamma^i_{jk}\frac{ds^j}{dt}\frac{ds^k}{dt}=v
\frac{ds^i}{dt},
$$
where $s^i$ are local coordinates of $U$, and $t$ is a parameter,
which is called affine if $v=0$.

One can associate a second-order ODE to a projective structure by
picking a connection in the equivalence class, choosing local
coordinates $s^i=(x,y)$ say, and eliminating the parameter from
the geodesic
equations. The resulting equation determines the geodesics in terms of
the local coordinates, without the parameter. The equation is
as follows:

\be \label{eqnwithgammas}
\frac{d^2 y}{dx^2} = \Gamma^x_{yy} \Big( \frac{dy}{dx} \Big) ^3 + (2
\Gamma^x_{xy} - \Gamma^y_{yy}) \Big( \frac{dy}{dx} \Big)^2 + (\Gamma^x_{xx} -
2 \Gamma^y_{xy})\frac{dy}{dx} - \Gamma^y_{xx}.
\ee
A general projective structure is therefore defined by a second-order
ODE (\ref{secondorderode}).
In fact, two of the four functions $A_0,A_1,A_2,A_3$ can be eliminated by a 
coordinate transformation $(x,y) \rightarrow (\hat{x}(x,y),\hat{y}(x,y))$ 
which introduces two arbitrary functions.

On $TU$, the horizontal lifts of $\p / \p s^i$ are defined by
$$
S_i = \frac{\p}{\p s^i} - \Gamma^j_{ik} v^k \frac{\p}{\p v^j},
$$
where $v^i$, $i=1,2$ are the fibre coordinates of $TU$.
The geodesics on $U$ lift to
integral curves of the following spray on $TU$:
\begin{equation} \label{projspray}
\cS=v^i S_i=v^i \frac{\partial}{\partial s^i} - \Gamma^i_{jk}v^j v^k
\frac{\partial}{\partial v^i},
\end{equation}
Now $\cS$ is homogeneous of degree 1 in the $v^i$, so it projects to a section of 
a one dimensional distribution on $PTU$.
$PTU$ is the quotient of $TU-\{
0\text{-section} \}$
by the vector field $v^i \frac{\p}{\p v^i}$. If $\lambda$ is a standard
coordinate on one patch of the $\RP$ factor
\footnote{By standard coordinates $\lambda, \tilde{\lambda}$ on $\RP$
or $\CP$, we mean the usual coordinates $v^1/v^0$ and $v^0/v^1$, where $v^0,v^1$ are
homogeneous coordinates.},
then the spray has the form
\begin{equation} \label{nonhomspray}
\cS = \partial_x + \lambda \partial_y + (A_0(x,y) + \lambda A_1(x,y) +
\lambda^2 A_2(x,y) + \lambda^3 A_3(x,y)) \partial_{\lambda}.
\end{equation}
There is a unique curve in any direction through a point in $U$, which
is why the curves can be lifted to a foliation of the projective
tangent bundle $U \times \RP$.

To obtain (\ref{connectiondifference}) we argue as follows. If
$\tilde{\cS}$ is the spray corresponding to a different connection
$\tilde{\Gamma}$, then $\Gamma$ and $\tilde{\Gamma}$ are in the same
projective class if $\cS$ and $\tilde{\cS}$ push down to the same spray
on $PTU$. This gives 
$$
\cS - \tilde{\cS} \propto v^i \frac{\p}{\p v^i},
$$ 
from which (\ref{connectiondifference}) follows, using the fact that
the connections are torsion-free (i.e. symmetric in their lower indices).

\subsection{The twistor space of a projective structure} \label{projtwistor}

Now suppose we have a \emph{holomorphic} projective structure on a
local patch of $\mathbb{C}^2$, which we still denote $U$. All of the
above is still valid, with
real coordinates replaced by complex ones. The functions
$\Gamma^i_{jk}$ are now required to be holomorphic functions of the
coordinates. Given a real-analytic projective structure, one can
complexify by analytic continuation to obtain a holomorphic projective
structure that will come equipped with a \emph{reality structure} (see
below).

The space $PTU$ is obtained from $TU$ on quotienting by $\mu^i \frac{\p}{\p
\mu^i}$, which defines a tautological line bundle
$\OO(-1)$ over $PTU$. 

As the $S_i$ are weight zero in the $\mu^i$ coordinates,
they push down to
vector fields on $PTU$, giving a two-dimensional distribution
$\mathcal{S}$. Since $\cS$ is weight one in the $\mu^i$, one must
divide by a homogeneous polynomial of degree one in the $\mu^i$ to
get something that pushes down to a vector field on $PTU$. The
resulting vector field will have a singularity at a single point on
each fibre, where the degree one polynomial vanishes. Different
choices of polynomial will result in different vector fields on $PTU$,
but they will always be in the same direction. In other words, $\cS$
defines a one dimensional distribution which we shall call  
$D_\cS$. Restricting to a
$\CP$ fibre, $D_\cS$ defines a line bundle over $\CP$. A section of this
line bundle corresponds to a vector field in $D_\cS$, i.e. a choice of
polynomial as described above, and has a pole at a single
point. Therefore by the classification of holomorphic line bundles
over $\CP$, it must be $\OO(-1)$ \symbolfootnote[2]{Coordinatize $\CP$
using two patches, $\mathcal{U}_0$ with coordinate $\lambda \in
\mathbb{C}$, and $\mathcal{U}_1$ with coordinate $\eta
\in \mathbb{C}$, and transition function $\lambda = 1 / \eta$. The
holomorphic line bundle $\OO(n)$ over $\CP$ is defined
by the transition function $a=\lambda^n b$ where $a \in \mathbb{C}$ 
is the fibre
coordinate over $\mathcal{U}_0$ and $b(\eta) \in \mathbb{C}$ is the
fibre coordinate over
$\mathcal{U}_1$. The
Birkhoff-Grothendieck theorem states that any holomorphic
line bundle over $\CP$ is $\OO(n)$ for some $n$. A global section of
$\OO(n)$ has $|n|$ zeroes or poles, for $n$ positive or negative respectively.}.

Restricting to a $\CP$ fibre,
one obtains the following exact sequence of vector bundles over $\CP$
\begin{equation} \label{projexactsequence}
0 \rightarrow \OO(-1) \rightarrow \OO \oplus \OO \rightarrow
\mathcal{S} / D_\cS \rightarrow 0,
\end{equation}
where the first bundle is $D_\cS$, the second is
$\mathcal{S}$, and the last is the quotient. In fact, the quotient
is $\OO(1)$, for the following reason. Consider for instance the push
down of $S_0$ to
$PTU$. This defines a subbundle of $\mathcal{S}$ that is different to
$\cS$ everywhere except at a single point, the image of $\mu^1 =
0$. Hence it determines a section of $\mathcal{S} / \cS$ which
vanishes at a single point. Therefore, again using the classification
of holomorphic line bundles over $\CP$, we have $\mathcal{S} / D_\cS
\cong \OO(1)$.

The twistor space $Z$ is the two dimensional quotient of $PTU$ by
$D_\cS$. A
point $u \in U$ corresponds to a twistor line $\hat{u} \subset Z$
corresponding to
all the geodesics through $u$.
The normal bundle of an embedded $\hat{u}=\CP \subset Z$ is
given by the quotient bundle in the above sequence, i.e. $\OO(1)$.
This is summarized by a double fibration picture:
$$
\begin{array}{lcccr}
&& U \times \CP&&\\
& \ \ \ \swarrow&&\searrow &\\
&U&& \ Z&
\end{array}
$$
The left arrow denotes projection to $U$, and the right arrow denotes
the quotient by $D_\cS$. 

The converse is also valid:
\begin{theorem}{\cite{hitchin,lebrun}}
\label{hitchin_lebrun}
There is a 1-1 correspondence between local two dimensional
holomorphic projective structures and complex surfaces containing an
embedded $\CP$ with normal bundle $\OO(1)$.
\end{theorem}
A vector $V \in T_u U$ corresponds to a global section of the normal
bundle $\OO(1)$ of $\hat{u}$. Such a section vanishes at a single
point $p \in Z$. The geodesic of the projective structure through this
direction is given by points in $U$ corresponding to twistor lines in $Z$
that intersect $\hat{u}$ at $p$. That there is a
one-parameter family of such lines can be shown by blowing up $Z$ at
the vanishing point and using Kodaira theory, see \cite{hitchin}. 

\subsection{Flatness of projective structures} \label{flatness}
A projective structure is said to be \emph{flat} if the corresponding
second order ODE (\ref{secondorderode}) can be transformed to the trivial ODE
\begin{equation} \label{trivialODE}
\frac{d^2 y}{dx^2} = 0
\end{equation}
by coordinate transformation $(x,y) \rightarrow (\hat{x}(x,y),
\hat{y}(x,y))$. The terminology comes from the fact that given
any second order ODE one can construct a Cartan connection on a
certain $G$-structure \cite{BGH}, and when this
connection is flat the equation can be transformed to the trivial ODE
(\ref{trivialODE}). It turns out that a second order ODE must be of
the form (\ref{secondorderode}) to be flat, and in addition the
functions $A_0, A_1, A_2, A_3$ must satisfy some PDEs. Defining 
$$
F(x,y,\lambda)=A_0(x,y) + \lambda A_1(x,y) + \lambda^2 A_2(x,y) + 
\lambda^3 A_3(x,y),
$$
the following must hold \cite{BGH}:

\be
\label{flatprojstructure}
\frac{\d^2}{\d x^2}F_{11}-4\frac{\d}{\d x}F_{01}-F_1\frac{\d}{\d x}F_{11}
+4F_1F_{01}-3F_0F_{11}+6F_{00}=0,
\ee
where 
\[
F_0=\frac{\p F}{\p y}, \qquad F_1=\frac{\p F}{\p \lambda}, \qquad
\frac{\d}{\d x}=\frac{\p }{\p x}+\lambda \frac{\p }{\p y}
+F\frac{\p}{\p \lambda}
\]
This is a set of PDEs for the functions $A_0, A_1, A_2, A_3$.

\subsection{Reality conditions for projective structures}
A reality structure for $Z$ is an anti-holomorphic involution that
leaves invariant a two real parameter family of twistor lines, and
fixes an equator of each line. Given a
line in this real family, all the sections pointing to nearby lines in
the real family have a zero at some point, and the union of these
points gives an equator of the line; this equator must be fixed by the
reality structure. The real family of twistor lines then corresponds
to a real manifold $U$ with a projective structure. 

In this paper all holomorphic projective structure have reality
structures since they occur as complexifications of real projective
structures.

\section{Null Killing Vectors and Twistor Space} \label{main}
\setcounter{equation}{0}
\subsection{The twistor space of an ASD conformal structure}

In the following and for the rest of the paper, $\et_{AA'}$ denote the
horizontal lifts of $\e_{AA'}$ to $S'$, or their push-down to $PS'$.
Abstractly, the integral curves of these horizontal
lifts define parallely transported primed spinors using the connection on
$S'$ (see Appendix \ref{appendix}).

We can abstractly define the two-dimensional \emph{twistor
  distribution} on $S'$ as follows. A
point $s \in S'$ is determined by a primed spinor $\pi$
at a point $x \in M$. The null vectors $\pi \otimes \mu$ for all
unprimed spinors $\mu$ span an $\alpha$-plane at $x$. Define the
twistor distribution at
$s$ to be the subspace of horizontal vectors at $s$ whose push-down to
the base lies in this $\alpha$-plane. 

Concretely, the twistor distribution is spanned by vectors $L_A$
$(A=0,1)$ on $S'$, defined with a choice of tetrad by 
\begin{equation} \label{twistordistribution}
L_A = \pi^{A'}
\et_{AA'} = \pi^{A'} \big( \e_{AA'} - \Gamma_{AA'B'}^{\ \ \ \ \
  \ C'} \pi^{B'}\frac{\p}{\p \pi^{C'}} \big),
\end{equation} 
where $\pi^{A'}$ are the local
coordinates on the fibres of $S'$. In the Appendix it is shown that
the twistor distribution is integrable for ASD conformal structures,
which is a seminal result of Penrose \cite{penrose}. In other
words, given a neutral ASD conformal structure $[g]$, each self-dual two plane
at a point is tangent to a unique \as \ through that
point, which is the push down of a leaf of the twistor distribution.
In the holomorphic case, the
space of leaves of the twistor distribution (locally, over a
suitably convex region of the base), is a three dimensional complex
manifold $\PT$ called the twistor
space \cite{penrose, hitchin}. 

The double fibration picture is very similar to the projective
structure case discussed in Section \ref{projtwistor}. The projective
primed spin bundle $PS'$ is the quotient of $S'$ by the vector field
$\pi^{A'} \frac{\p}{\p \pi^{A'}}$. $PS'$ is fundamental in the fibration
picture, as each \as in $M$ has a unique lift, in the same
way that each geodesic of a projective structure has a unique lift to
the projective tangent bundle. The horizontal
vectors $\et_{AA'}$ are weight zero in the $\pi^{A'}$ coordinates, so
push down to vector fields on $PS'$, giving a four-dimensional
distribution $\Xi$ on $PS'$. The $L_A$ vectors
(\ref{twistordistribution}) are weight one, so together define a two
dimensional subdistribution $\cL$ of $\Xi$, which restricts to $\OO(-1)
\oplus \OO(-1)$ on a $\CP$ fibre; we also refer to this as the twistor
distribution. Over a $\CP$ fibre, there is an
exact sequence
\begin{equation} \label{twistorsequence}
0 \rightarrow \OO(-1) \oplus \OO(-1) \rightarrow \OO \otimes
\mathbb{C}^4 \rightarrow \Xi / \cL  \rightarrow 0.
\end{equation}
The first term is $\cL$, the second is $\Xi$.
As in the projective structure case, one can show that $\Xi /
\cL$ is $\OO(1) \oplus \OO(1)$. The twistor space $\PT$ is the
quotient of $PS'$ by $\cL$. The image of a $\CP$ fibre over $x \in M$ is an
embedded $\CP \in \PT$, and has normal bundle $\OO(1) \oplus \OO(1)$,
the quotient bundle in (\ref{twistorsequence}). It corresponds to all
the \as s through $x$. 

The twistor correspondence is summarized by the double fibration:
$$
\begin{array}{lcccr}
&& PS'&&\\
& \ \ \ \swarrow&&\searrow &\\
&M&& \PT&
\end{array}
$$
Here the left arrow denotes projection to $M$, and the right arrow
denotes the quotient by $\cL$. 

Again, there is a converse:
\begin{thm}{(Penrose \cite{penrose})}
There is a 1-1 correspondence between local four dimensional holomorphic ASD
conformal structures $(M,[g])$ and three dimensional complex manifolds
$\PT$ with an embedded $\CP \subset \PT$, with normal bundle $\OO(1)
\oplus \OO(1)$.
\end{thm}
The essential fact is that an embedded $\CP$ with the above normal bundle belongs
to a family of embedded $\CP$s parameterized by a complex 4-manifold $M$.
Vectors at $x \in M$ correspond to sections of the normal bundle of
$\hat{x}$, and null vectors are given by sections with a zero. This
defines a conformal structure, because a global section of $\OO(1) \oplus
\OO(1)$ is given by $(a \lambda + b, c \lambda + d)$ for affine
coordinate $\lambda \in \C$, $(a,b,c,d) \in \C^4$, and this can only
be $(0,0)$ when $ad-bc=0$,
which is a quadratic condition. In this case there is a zero at a single point.
The conformal structure is anti-self-dual, with \as s defined by
families of twistor lines through a fixed point in $\PT$.

In this picture, the \as s are obtained as follows. Let $\hat{x}
\subset \PT$ be the twistor line corresponding to a point $x \in
M$. Let $V \in T_x M$ be a null vector. We want to show that $V$ lies in a
unique \as \ through $x$. The corresponding section of the normal bundle of
$\hat{x}$ has a zero at some point $p \in \PT$ because $V$ is null. The \as
\ corresponds to all the twistor lines that intersect $\hat{x}$ at $p$.
There is a two-parameter family of such lines. It is easy to see that
there is a two-parameter family of sections that vanish at $p$. To
show that these are tangent to a two-parameter family of lines one
must blow-up $\PT$ at $p$ and use Kodaira theory; see
\cite{hitchin} for details.

\subsubsection{Reality conditions for split signature}
In order to obtain a real split signature metric from a twistor space,
we must be able to distinguish a four \emph{real} parameter family of
twistor lines, whose parameter space will be the four real dimensional
manifold. In
addition we require that given a line in this real family, the sections of
the normal bundle that point to others in the family inherit a
split signature conformal structure. As described above, a section of
$\OO(1) \oplus \OO(1)$ is defined by four complex numbers $(a,b,c,d)$,
with a quadratic form defined by $ad-bc$. If we restrict $(a,b,c,d)$
to be real we obtain a real split signature conformal structure. The
sections tangent to the real family are of this type.

The zero of such sections occurs for real $\lambda$, that is, on an
equator of $\CP$. The conformal structure is thus invariant under an
anti-holomorphic involution of the $\CP$ that has the equator fixed. A
\emph{split signature real structure} on $\PT$ is an anti-holomorphic
involution that leaves invariant a four real parameter family of
twistor lines, and when restricted to one of these fixes an equator. 

Not all holomorphic metrics have real structures, but all the
holomorphic metrics in this paper have obvious real
`slices' because they are complexifications of real metrics, obtained
by letting the real coordinates be complex.

\subsection{Lift of $K$ to $PS'$}

Now given a null conformal Killing vector for an ASD conformal structure,
the fact that $M$ is foliated by \as s  (Lemma \ref{kvproperties}) is
not very illuminating, since they must already exist by anti-self-duality. The
foliation by \bs s is more interesting, since these do not exist generically.

In this section we will prove that in the analytic case, the space of
\bs s inherits a natural projective structure. We then explain
how this arises geometrically, due to the presence of \as s ensured by
anti-self-duality.

Let $K$ be a null conformal Killing vector for $(M,[g])$. We assume
$K$ is without
fixed points, which can always be arranged by restricting $M$ to a
suitable open set.

Since $K$ preserves the conformal structure, the corresponding diffeomorphism  maps \as s to
\as s, and hence it induces a vector field $\cK$ on $\PT$.
We now translate this fact into a statement on the projective primed
spin bundle $PS'$. Each \as \ has a
unique lift and these lifts foliate $PS'$. The following proposition
shows how to lift $K$ to $PS'$, giving a vector field that is
Lie-derived along the lifts of the \as s.

\begin{proposition} \label{liftofK}
Let $K=K^{AA'}\e_{AA'}$ be a conformal Killing vector for an ASD
metric $g$. Define a vector field $\tK$ on $S'$ by 
\begin{equation} \label{Klift}
\tK:=K^{AA'}\et_{AA'}+\pi_{A'}\phi^{A'B'}\frac{\p}{\p
  \pi^{B'}}+\frac{1}{2} \eta \pi^{A'} \frac{\p}{\p \pi^{A'}}.
\end{equation}
Then this satisfies 
\begin{equation} \label{lielift}
[\tK,L_A] = (K^{BB'} \Gamma_{BB'A}^{\ \ \ \ \ \ D}-\psi_A^{\ \
  D})L_D + \frac{3}{4} (\e_{AB'} \eta) \pi^{B'} \pi^{C'} \frac{\p}{\p
 \pi^{C'}}.
\end{equation}
\end{proposition}

\emph{Proof.} See Appendix. $\square$
\vspace{10pt}

\textbf{Remark.} Since $\tK$ is weight zero in the $\pi^{A'}$
coordinates, it defines a vector field on $PS'$, which we will also
refer to as $\tK$ by abuse of notation. The last term on the
right hand side of (\ref{lielift}) is proportional to the Euler vector
field, so does not contribute to $\tK$ on $PS'$. Hence (\ref{lielift})
shows that $\tK$ commutes with the twistor distribution $\cL$ on
$PS'$. The vector field $\cK$ on $\PT$ is the push-forward of $\tK$ to
$\PT$, which is well defined because $\tK$ is Lie-derived along $\cL$.

\subsection{Projective structure from a quotient} \label{projquotient}

In this section we assume
that $[g]$ is analytic, so we can complexify by analytic
continuation. Thus we are now working on a local patch of
$\mathbb{C}^4$, with a holomorphic conformal structure.
We assume that we have restricted to a
suitable open set on the base so that all the spaces of leaves involved
are non-singular complex manifolds. 

As in Section \ref{NullKVs}, write
$K=\iota^A o^{A'} \e_{AA'}$, where now $\e_{AA'}$ is a holomorphic
tetrad and $\iota^A$, $o^{A'}$ are complex spinor fields that vary
holomorphically.

When $K$ is null,
it is easy to see that $\cK$, the induced vector field on twistor
space $\PT$, will vanish on a hypersurface $\mathcal{H}$ in $\PT$,
because $K$ fixes a two-parameter family of \as s, which are those to
which it is tangent. These are the `special' \as s of Lemma
\ref{kvproperties}. We now explain how this is seen from the lift
$\tK$ to $PS'$.

On $S'$,
$\tK$ from Proposition {\ref{liftofK}} is given by:
$$
\tK = \iota^A o^{A'} \et_{AA'} + \pi_{A'}\phi^{A'B'}\frac{\p}{\p
  \pi^{B'}} + \frac{1}{2} \eta \pi^{A'} \frac{\p}{\p \pi^{A'}}.
$$
Now when $\pi^{A'} \propto o^{A'}$, one has $\pi_{A'} \phi^{A'B'} \propto
o^{B'}$ from (\ref{alg1}), so the second term on the RHS is
proportional to the Euler
vector field $\Upsilon = \pi^{A'} \frac{\p}{\p \pi^{A'}}$. The last term is
everywhere proportional to the Euler vector field. To go from $S'$ to
$PS'$ one quotients $S'-\{0 \text{-section} \}$ by the integral curves of
$\Upsilon$. So we have shown that on the section $[\pi^{A'}]=[o^{A'}]$
of $PS'$, $\tK$ is the push down of $\iota^A o^{A'} \et_{AA'}$ only. But
this is in $\cL$, so $\tK$ pushes down to the zero vector under the
quotient of $PS'$ by $\cL$.

So there is a (complex) hypersurface in $PS'$, defined by
the section $[\pi^{A'}] = [o^{A'}]$, on which
$\tK$ lies in the twistor distribution. One can also define this
hypersurface as the image in $PS'$ of the hypersurface $\pi.o = 0$ in
$S'$, under the quotient by $\Upsilon$. We will refer to this
hypersurface as $H$. It is easy to see by pushing
down to the base that $\tK$ is linearly independent
of the twistor distribution everywhere else on $PS'$. 

Define a vector field 
\[
V=\iota^A L_A = \iota^A \pi^{A'} \et_{AA'}\]
 on
$S'$. This is weight one in the $\pi^{A'}$ coordinates, so gives a one dimensional
distribution on $PS'$ which restricts to $\OO(-1)$ on fibres. Together
with $\text{span} \{ \tK \}$, we get a two dimensional distribution on $PS' - H$. On
$H$, the distribution drops its rank from two to one. 

The two dimensional distribution defined by $\{V, \tilde{K}\}$
on $PS'-H$ pushes down to the
$\beta$-plane distribution $\D_{\beta}$ on the base.
\begin{lemma} \label{commutingdistribution}
The two dimensional distribution on $PS'-H$ determined by
$\{V, \tilde{K}\}$ is integrable.
\end{lemma}
\emph{Proof.} We work on $S'$ for convenience, and push down to $PS'$ at the
end. The distribution $\text{span} \{ \tK, V \}$ on $S'$ is two
dimensional on $S'$ when $\pi^{A'} o_{A'} \neq 0$. Multiples of the
Euler field $\Upsilon$
are therefore irrelevant.
\begin{eqnarray*}
[V,\tK] &=& [\tK, \iota^C L_C ] \\
        &=& \iota^C [\tK,L_C] + \tK(\iota^B) L_B \\
        &=& \iota^C ((K^{BB'} \Gamma_{BB'C}^{\ \ \ \ \ \ D}-\psi_C^{\ \
  D})L_D + \frac{3}{4} (\e_{CB'} \eta) \pi^{B'} \pi^{C'} \frac{\p}{\p
 \pi^{C'}}) \\ && \ \ + K^{BB'} \e_{BB'}(\iota^C) L_C \\
        &=& (K^{BB'} \nabla_{BB'} \iota^C - \iota^D \psi_D^{\ \ C})
        L_C + \# \Upsilon \\
        &=& (\iota^B o^{B'} \nabla_{BB'} \iota^C - \iota^D \psi_D^{\ \ C})
        L_C + \# \Upsilon.
\end{eqnarray*}

From (\ref{alg1}) we have $\iota^D \psi_D^{\ C} \propto \iota^C$, and
from ($\ref{GSF1}$) we have $\iota^B o^{B'} \nabla_{BB'} \iota^C \propto
\iota^C$. Hence the RHS is proportional to $V$, ignoring the
irrelevant Euler vector field part. $\square$

Next we will show that it is possible to continue this distribution
 over the hypersurface $H$ so it is rank two on the whole of $PS'$, and that the
resulting distribution commutes on the hypersurface. It will then be
possible to quotient $PS'$ by the leaves of this distribution.

\begin{lemma} \label{hypersurface}
There is a two-dimensional integrable distribution $\D$ over $PS'$,
which on $PS'-H$
is determined by $\{ \tK, V\}$. Let $\varrho$ be the
projection $PS' \rightarrow M$. Then for
every $p \in PS'$, we have $\varrho_* (\D \mid_p ) = \D_{\beta}$.
\end{lemma}

\textbf{Remark.} Intuitively one can think
of $\D$ as a lift of the \bs s to $PS'$, where each \bs \ has a $\CP$ of lifts.

\vspace{10pt}

\emph{Proof.} Choose a spinor $\iota^{A'}$ satisfying $o^{A'}
\iota_{A'} = 1$. Define the following (singular) vector field on $S'$:
\begin{equation} \label{Wvector}
W = \frac{1}{\pi^{C'} o_{C'}}(V-(\pi^{D'}\iota_{D'}) \tK).
\end{equation}
This is weight zero in the $\pi^{A'}$, so defines a vector field on
$PS'$ by push-forward, which we shall also call $W$. We will now show
that $W$ is well defined even over $H \subset PS'$, despite the $1 /
(\pi^{C'}o_{C'})$ factor in (\ref{Wvector}). 

Without loss of generality, choose a tetrad such that
$$K = \iota^A o^{A'} \e_{AA'} =\e_{00'}.$$ That is, $\iota^A = (1 ,
0)$, $o^{A'}
= (1 , 0)$. Define $\lambda = \pi^{1'}/\pi^{0'}$ to be the
coordinate on the $\pi^{0'} \neq 0$ patch of $\CP$, and extend this to
a patch of $PS'$; we call the patch $\mathcal{U}$. Then $H$ lies
entirely within $\mathcal{U}$ at $\lambda = 0$. We have the following
expression
for $\tK$, obtained by `projectivizing' (\ref{Klift}):
\begin{eqnarray*}
\tK &=& \et_{00'} + (\phi_{0'}^{\ \ 1'} + \lambda (\phi_{1'}^{\ \ 1'} -
\phi_{0'}^{\ \ 0'}) + \lambda^2 \phi_{1'}^{\ \ 0'}) \frac{\p}{\p
  \lambda} \\
    &=& \et_{00'} + (\lambda (\phi_{1'}^{\ \ 1'} -
\phi_{0'}^{\ \ 0'}) + \lambda^2 \phi_{1'}^{\ \ 0'}) \frac{\p}{\p
  \lambda},
\end{eqnarray*}
where $\phi_{0'}^{\ \ 1'} = 0$ due to (\ref{alg2}). 

In the above conventions, we have 
$V = \pi^{A'} \et_{0A'}$. On $\mathcal{U} \subset PS'$, the push
forward of $\frac{1}{\pi^{C'} o_{C'}}V$ is
$$
\frac{1}{\lambda}\et_{00'} + \et_{01'},
$$
which is singular at $H$, corresponding to $\lambda=0$. Choosing
$\iota^{A'}=(0,-1)$, we then obtain the
following expression for $W$ on $\mathcal{U}$:
$$
W = \frac{1}{\lambda} \et_{00'} + \et_{01'} - \frac{1}{\lambda} \tK = \et_{01'} - ((\phi_{1'}^{\ \ 1'} -
\phi_{0'}^{\ \ 0'}) + \lambda \phi_{1'}^{\ \ 0'})) \frac{\p}{\p
  \lambda}.
$$
This is a \emph{non-singular} vector field on $\mathcal{U}$.
By construction, away from $H$ this lies in $\text{span}\{
\tK, \tilde{V} \}$. Define $\D$ on $\mathcal{U}$ to be $\text{span} \{ \tK, W
\}$. This is clearly non-degenerate everywhere on $\mathcal{U}$. Note that
$W$ is also well defined over the other patch (i.e. at $\lambda =
\infty$) so we can define $\D$ as $\text{span} \{ \tK, W
\}$ over the whole of $PS'$.

We now want to show that $\D$ is integrable over $H$. We know (Lemma
\ref{commutingdistribution}) that $\D$ is integrable away from
$H$. Therefore on $\mathcal{U}$ we have 
$$
[\tK,W] = f \tK + g W + Y,
$$
where $f,g$ are holomorphic functions on $\mathcal{U}$ and $Y$ is a holomorphic
vector field vanishing on $\mathcal{U} - H$. But such a vector field must
vanish, otherwise it is not even continuous, so is certainly not holomorphic.

The last part of the lemma is obvious, just from inspecting the
coordinate expressions of $\tK$, $W$. $\square$

\vspace{10pt}

We now have a three dimensional integrable distribution $\cL \cup \D$. It is three
dimensional because at each point $\cL$ and $\D$ have a direction in
common, which is the one-dimensional distribution defined on $PS'$ by
the push-forward of $V$ on $S'$. From Lemma \ref{hypersurface}, $\D$
is an integrable subdistribution. Note that $\D$ consists of a $\CP$
of lifts of each \bs
\ in the base. If we pick a suitably convex set on the base so that
the space of \bs s $U$ intersecting it is a Hausdorff complex manifold,
then the quotient $PS' / \D$ will also be a
Hausdorff complex manifold. A point in this quotient is a point in $U$
together with a choice of lift. 

In fact we can canonically identify $PS' / \D$
with $PTU$, the projective tangent bundle of $U$, as follows. Using
the conventions of Lemma \ref{hypersurface}, the tangent planes to the
\bs s in the base are spanned at each point by $\e_{00'}$,
$\e_{01'}$. Now $L_1$ has the form $\e_{10'} + \lambda \e_{11'} +
(\ldots) \partial_{\lambda}$, so at each point in the fibre above a
point $x \in M$, $L_1$ pushes down to a different null direction
transverse to the $\beta$-plane at $x$. Now suppose we take a lift of
a \bs \ $\Pi$, i.e. a leaf of $\D$ that projects down to $\Pi$. Push down
$L_1$ at each point over this lift. This
will give a vector field $\Theta = \e_{10'} +
\lambda \e_{11'}$ over $\Pi$, where $\lambda$ is now a function on the
$M$.

We want to show that this determines a
projective vector at the point $s \in U$ corresponding to $S$. This
means we require $[\e_{00'},\Theta] \propto \Theta \ \text{mod} \{
\e_{00'}, \e_{01'} \}$, $[\e_{01'},\Theta] \propto \Theta \ \text{mod} \{
\e_{00'}, \e_{01'} \}$. But it is easy to show that this is satisfied,
using the fact that the distribution spanned by $\tK, W, L_1$
commutes. Hence to determine the projective vector corresponding to a
leaf of $\D$, just choose a point on the leaf and push down
$L_1$. Because of the above considerations, this direction will be
independent of the choice of point on the leaf.

\vspace{15pt}

\emph{Proof of Theorem \ref{maintheorem}.}  
Define $Z$ as the quotient of $PS'$ by $\cL \cup
\D$. Equivalently, this is the quotient of $\PT$ by a one-dimensional
distribution which on $\PT - \mathcal{H}$ is 
$\text{span} \{ \cK \}$. The image of a $\CP$ fibre of $PS'$ under the
quotient is a twistor line in $Z$.

On a $\CP$ fibre, the horizontal part of $\D$ defines a subbundle $\OO
\otimes \mathbb{C}^2$ of the horizontal distribution $\Xi = \OO
\otimes \mathbb{C}^4$, corresponding to the horizontal parts of the
vectors $\tK$ and $W$. Choosing a spinor $o^A$ such that $\iota^A o_A
= 1$, we can form the vector field $o^A L_A$ on $S'$, which pushes
down to a horizontal distribution on $PS'$ that is always linearly independent of
$\D$. Since the $L_A$ are weight one, this is $\OO(-1)$ when
restricted to a $\CP$ fibre. Because $\cL \cup \D$ is integrable
(Lemma \ref{hypersurface}), this
distribution determines a one dimensional distribution  $D_\Theta$ 
on $PTU = PS' / \D$. The spray $\cS$ of a projective structure
is a section of $\D_\Theta\otimes \OO(1)$ where here $\OO(1)$ is dual to 
the tautological line bundle over $PTU$.
The situation is described by the following commuting diagram:
$$
\begin{array}{lcccccccr}
0&\rightarrow& \OO(-1)\oplus \OO(-1)& \rightarrow & \OO \otimes \mathbb{C}^4 &
\rightarrow & \OO(1)\oplus \OO(1) & \rightarrow & 0 \\
&&\downarrow \scriptstyle V&&\downarrow \scriptstyle \D&&\downarrow
&& \\
0&\rightarrow & \OO(-1)& \rightarrow & \OO \otimes \mathbb{C}^2 &
\rightarrow & \OO(1) & \rightarrow & 0
\end{array}
$$
where these are bundles over a $\CP$ fibre of $PS'$. The vector field
$o^A L_A$ on $S'$ constructed above corresponds to
the $\OO(-1)$ in the bottom row after quotienting by $V$, and gives
the projective structure spray. The bottom row is the sequence
(\ref{projexactsequence}) on $PTU = PS' / \D$, where $U$ is the space
of \bs s in $M$. Thus there is a projective structure on $U$. $\square$

\vspace{10pt}
\textbf{Remark.} The real space of \bs s has a system of curves that
comes from the quotienting operations described above but with real
spaces instead of complex. These real curves are described by the
holomorphic projective structure with a reality structure.
\vspace{10pt}
\begin{figure} \label{fibrediagram}
\centering
\includegraphics[height=1.5in]{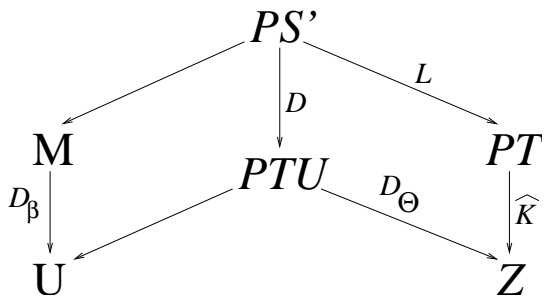}
\caption{Relationship between foliation spaces.}
\end{figure}

Figure 1 illustrates the situation. Here $p$ and $q$ are the obvious
projections. $\mathcal{D}_{\beta}$ represents the \bs \
distribution on $M$. The $\widehat{\cK}$ labelling the map from $\PT$ to $Z$
requires some of explanation. The vector field $\tK$ over $PS'$
commutes with the twistor distribution (Lemma \ref{liftofK}), so
determines a vector field $\cK$ on $\PT$. This vector field vanishes
on a hypersurface $\mathcal{H} \subset\PT$, corresponding to the \as s to
which $K$ is tangent; these are the \as s appearing in Lemma
\ref{kvproperties}. Now $\cK$ on $\PT$ only depends on $\tK$ modulo
$\mathcal{L}$. Lemma \ref{hypersurface} shows that we can multiply
$\tK$ modulo $\mathcal{L}$ by a meromorphic function ($1/ \lambda$)
and obtain a vector field $W$ commuting with the twistor
distribution. This means that there is a one-dimensional
distribution $\widehat{\cK}$ over the whole of $\PT$ that never
degenerates, and which
agrees with $\text{span} \ \{ \cK \}$ on $\PT - \mathcal{H}$. The quotient of
$\PT$ by this distribution gives $Z$, as illustrated in the diagram.

One can rephrase this in terms of divisor line bundles. That is, there
is a holomorphic line bundle $\mathcal{E}$ over $\PT$ defined by the
property that it has a meromorphic section $\zeta$ with a pole of
order one on $\mathcal{H}$. Then $\zeta \otimes \cK$ defines a
non-vanishing section of $\mathcal{E} \otimes T \PT$. This is
equivalent to the one dimensional distribution $\widehat{\cK}$ over
$\PT$ described
above. To obtain the distribution one simply finds trivializations of
$\mathcal{E}$ and $T \PT$ over a patch, and expresses $\zeta$ in this
trivialization. Its direction will be independent of the
trivialization of $\mathcal{E}$, and defines the distribution over the
patch.

\vspace{10pt}
\subsection{Relationship of the twistor spaces}
\begin{figure}
\caption{The $\alpha$ and $\beta$ surfaces in $M$ intesect
along a trajectory of $K$ which is a null geodesic. This corresponds to
a point $\alpha$ lying on a surface $\beta$ in ${\cal PT}$. Points $p_1, p_2, p_3$ in $M$ correspond to projective lines in ${\cal PT}$.}
\label{intersections_fig}
\begin{center}
\includegraphics[width=11cm,height=8cm,angle=0]{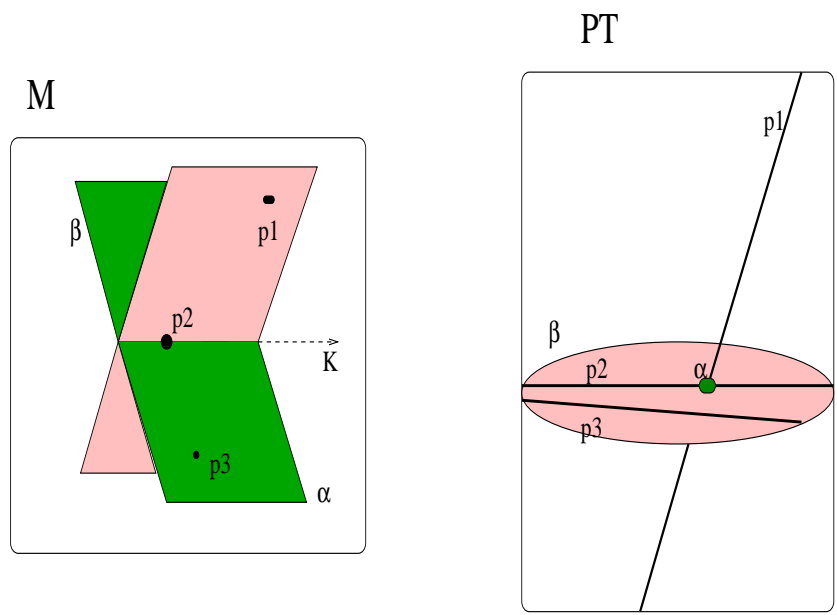}
\end{center}
\end{figure}

Here we discuss the relationship between the twistor
spaces without the foliation space picture. Incidence relation between various objects in $M$ and ${\cal PT}$ is represented by (Fig. \ref{intersections_fig}).

First one must
understand what a \bs \ corresponds to in $\PT$. The answer is a
two-parameter family of twistor lines, each of which intersects any
other at a single point. 
This is because all points on a \bs
\ are null separated. However, unlike the case of an \as, there is
not just a single point of intersection of the whole family. To
construct the family, pick a point on the \bs, say $x$. Then $\hat{x}$
is a twistor line in $\PT$. Now $\cK$ determines a section of the
normal bundle with a zero. Twistor lines intersecting $\hat{x}$ at
this zero are on the \bs, and correspond to those along the trajectory
of $K$ through $x$. In fact this is a null geodesic, since null Killing
vector fields have geodesic integral curves. Now pick another section
of the normal bundle with a zero at a different point, such that all
linear combinations of this
with the section determined by $\cK$ also have a zero. The resulting
two dimensional distribution in $M$ at $x$ is a $\beta$-plane. Doing
this for each $x \in M$ gives a $\beta$-plane distribution which is
integrable.

\begin{figure}
\label{Figg}
\centering
\includegraphics[height=8cm]{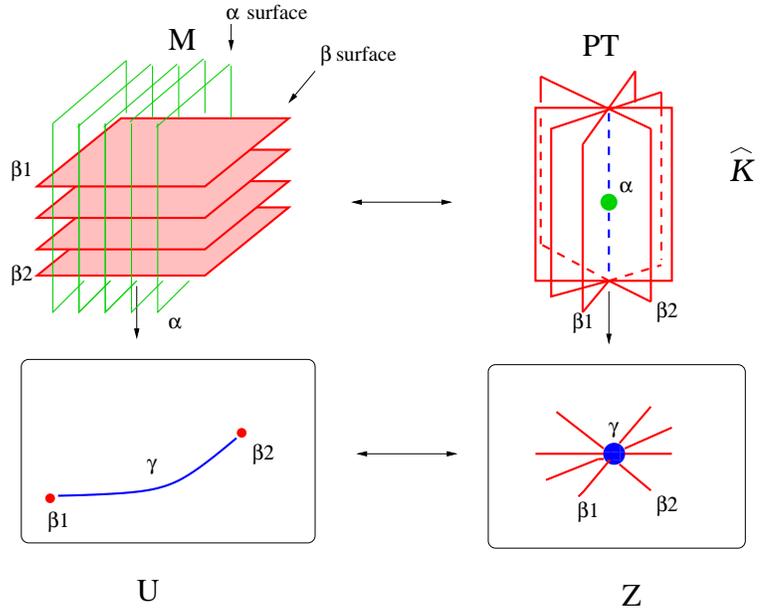}
\caption{Relationship between $M$, $U$, $\PT$ and $Z$.}
\end{figure}

The diagram (Fig. \ref{Figg})  illustrates the situation. 
In $M$, a one parameter family of
\bs \ is shown, each of which intersects a one parameter family of \as
s, also shown. The \bs s correspond to a projective structure geodesic
in $U$, shown at the bottom left.

The \bs s in $M$ correspond to surfaces in $\PT$, as discussed
above. These surfaces intersect at the dotted line, which corresponds
to the one parameter family of \as s in $M$. When we quotient $\PT$ by
$\cK$ to get $Z$, the surfaces become twistor lines in $Z$, and the
dotted line becomes a point at which the twistor lines intersect; this
is shown on the bottom right. This family of twistor lines
intersecting at a point corresponds to the geodesic of the projective
structure.

\section{Local classification} \label{local}
\setcounter{equation}{0}
The second theorem stated in the Introduction gives a local expression
for any analytic neutral signature ASD conformal structure. We now
prove this theorem. In the proof we will often use the following shorthand
for coordinate transformations: $t \rightarrow F(t,x,y,z)$ means
define a new coordinate $\tilde{t} = F(t,x,y,z)$ and then relabel it
$t$ again. This avoids having to introduce new symbols for new
coordinates. We will denote partial derivatives by subscripts, for
example $F_z := \p_z F$.

\emph{Proof of Theorem \ref{localform}.} In what follows, we will use
coordinates $(x,y)$ for the two-dimensional space of \bs s $U$. We
will always work on a single patch of $PS'$, with $\lambda$ a standard
coordinate on one patch of the $\CP$ fibre. The projectivization of
(\ref{twistordistribution}) is
\begin{eqnarray}
L_0 &=& \e_{00'} + \lambda \e_{01'} + (f_0 + \lambda f_1 + \lambda^2
f_2 + \lambda^3 f_3) \p_{\lambda},  \label{projtwistordistribution1} \\
L_1 &=& \e_{01'} + \lambda \e_{11'} + (A_0 + \lambda A_1 + \lambda^2
A_2 + \lambda^3 A_3 ) \p_{\lambda},  \label{projtwistordistribution2}
\end{eqnarray}
where the $f_\alpha$ and $A_\alpha$ are functions on $M$ derived from primed
connection coefficients.
 
We can trivialize $PTU$ by first
choosing a two dimensional surface in $M$, transverse to the
\bs s, and trivializing $PS'$ over this, using the standard two patch
coordinates for $\CP$. Then define a trivialization over the rest of
$PS'$ by requiring constant coordinate on each leaf of
$\D$ (this will be a base dependent M\"obius transformation of any
other trivialization of $PS'$ using a standard two patch
trivialization of $\CP$, since any two standard
trivializations of $\CP$ are related by a M\"obius transformation). This gives
a trivialization $PTU \cong U \times \CP$. The special feature
of this particular trivialization is that $\tK$ and $W$ will have no
vertical terms, because it was defined by saying that the fibre
coordinate is constant along them.

We will use the conventions of Lemma \ref{hypersurface}, that is we
choose a tetrad with $K=\e_{00'}$, and the tangent planes to the \bs s
are spanned by $K$ and $\e_{01'}$. Now choose a coordinate system
$(t,x,y,z)$ such
that $K=\p_t$, and a conformal factor so that $K$ is pure Killing. Any
tetrad can then be written in these coordinates without any $t$
dependence. Then $[\e_{00'}, \e_{01'}]=0$ and we can in addition
choose the $z$ coordinate such that $\e_{01'} = \p_z$. Then we have
\begin{eqnarray*}
\tK &=& \partial_t, \\
L_0 &=& \p_t + \lambda \p_z + f(x,y,z,\lambda) \partial_{\lambda}.
\end{eqnarray*}
Note that $f$
does not depend on $t$ because it is composed from
connection coefficients, which do not depend on $t$ since it does not
occur in the metric. Also note that $\tK = \p_t$ because $L_0$, $L_1$
do not contain functions of $t$ so it commutes with both. As vector
fields on the base, $\p_x$ and $\p_y$ are transverse to the \bs s so
are coordinates on $U$.

Now we will alter the $\lambda$ coordinate, using a trivialization as
described above, so that $L_0$ has no $\p_{\lambda}$ terms. This
is achieved by a M\"obius transformation, $\lambda
\rightarrow (\beta + \delta \lam)/ (\alpha + \gamma \lam)$, where $\alpha, \beta, 
\gamma, \delta$ are functions on $M$. Now the
new $\lambda$ coordinate satisfies  $\tK(\lambda)=L_0
(\lambda)=0$. Therefore $\alpha, \ldots, \delta$ do not depend on $t$,
from the first of these. This gives the following general form:
\begin{eqnarray}
\tK &=& \p_t,  \label{Kform} \\
L_0 &=& \alpha \ \p_t + \beta  \ \p_z + \lam
(\gamma \ \p_t + \delta \ \p_z). \label{L0form}
\end{eqnarray}

Now from Theorem \ref{maintheorem}, we know that $L_1$ must define a
projective structure on $U$, the space of \bs s. In fact this can be 
seen directly using our coordinate choices \footnote{The following 
argument does not require analyticity, only smoothness. Consequently Theorem 
\ref{localform} will turn out to be valid for smooth conformal 
structures, not just analytic ones. The smooth generalisation of 
Theorem \ref{maintheorem} 
could perhaps be established using techniques 
introduced in \cite{lebrunmason}.}. Clearly $U$ has
coordinates $(x,y)$, since the \bs s are spanned by $(\p_t,\p_z)$. Also,
$\lam$ is a fibre coordinate on $PTU$, since $\D$ is defined by
constant $\lam$. Since $\{ L_{0}, L_{1} \}$ is an integrable 
distribution, one can find a non-zero function $f$ on $PS'$ such that 
$[L_{0},f L_{1}] \propto L_{0}$. We may therefore assume we have 
chosen an $L_{1}$ such that $[L_{0},L_{1}] 
\propto L_{0}$. It follows from (\ref{L0form}) that the coefficients 
in front of the $\p_{x}, \p_{y}, \p_{\lambda}$ terms in $L_{1}$ do not 
depend on $z$. Therefore $L_1$ must have the following form:
\begin{multline} \label{L1form}
L_1:=J_0(x,y)\partial_x + J_1(x,y) \p_{y}+ \lam (J_2(x,y) \partial_x 
+J_3(x,y) \p_{y}) \\ + (A_0(x,y) + \lam A_1(x,y) + \lam^2
A_2(x,y) + \lam^3 A_3(x,y))\p_{\lam} \\ +
(C(x,y,z) + \lam D(x,y,z))\p_t + (E(x,y,z) + \lam F(x,y,z)) \p_z,
\end{multline}
where $J_0J_3-J_1J_2 \neq 0$. One now observes that the $\p_{x}, \p_{y}, 
\p_{\lambda}$ terms precisely correspond to a projective structure 
spray on $PTU$. Since $\D$ is 
spanned  by $\p_t, \p_z$, the quotient of $L_1$
by $\D$ gives a projective structure.

To put the projective structure spray occuring in $(\ref{L1form})$ into 
the more standard form 
$(\ref{nonhomspray})$ (i.e. $J_0=J_3=1,J_1=J_2=0$) it is necessary 
to perform a M\"obius 
transformation of $\lambda$ depending on $(x,y)$. Since this does not 
depend on $t$ or $z$, the general form $(\ref{L0form})$ of $L_{0}$ is 
unchanged by this, and we can assume that the projective structure 
spray in $L_{1}$ is of the form $(\ref{nonhomspray})$, which we shall do 
from now on.

We have found a general form that any $\{ \tK,L_A \}$ can be put into. For
it to give an ASD conformal structure, the $L_A$ must
commute modulo $L_A$. Imposing this gives equations for the unknown
functions, which will lead us to the metrics appearing in Theorem 
\ref{localform}.

First, it is convenient to change coordinates yet again, because
together with conformal rescaling we can elimate
three of the four functions in $L_{0}$. We may assume 
$\delta \neq 0$ (if $\delta=0$ then $\beta \neq 0$, in which case
perform the coordinate change $\lambda \rightarrow 1 / \lambda$).

Now change 
coordinates by $({t},{x},{y},{z})\rightarrow(t + 
j(x,y,z),x,y,k(x,y,z))$, where $k_z\neq 0$.
%Then we obtain
%$$
%L_{0} = \alpha \p_{\tilde{t}} + \beta(\frac{\p k}{\p z} \p_{\tilde{z}} + 
%\frac{\p j}{\p z} \p_{\tilde{t}}) + \lambda (\gamma \p_{\tilde{t}}+ \delta 
%(\frac{\p k}{\p z} \p_{\tilde{z}} + 
%\frac{\p j}{\p z} \p_{\tilde{t}})).
%$$
%Here we regard $\frac{\p k}{\p z}$ and $\frac{\p j}{\p z}$ as 
%functions of $(\tilde{x},\tilde{y},\tilde{z})$.
%Now choose $j(x,y,z)$ such that the following equation is satisfied:
%$$
%\frac{\p j}{\p z} = - \frac{\gamma}{\delta}.
%$$
%The right hand side is not singular because $\delta \neq 0$. Then 
%we have 
%$$
%L_{0} = (\alpha - \frac{\beta \gamma}{\delta}) \p_{\tilde{t}} + \beta 
%\frac{\p k}{\p z} \p_{\tilde{z}} + \lambda \delta \frac{\p k}{\p z} \p_{\tilde{z}}
%$$
%Since the function $(\alpha - \frac{\beta \gamma}{\delta}) \neq 0$, we can divide 
%by it. Finally one can choose $k(x,y,z)$ to satisfy
%$$
%\frac{\p k}{\p z} = \frac{1}{\delta}(\alpha - \frac{\beta 
%\gamma}{\delta}),
%$$
%where the right hand side is non-vanishing.
A suitable choice of $j$ and
$k$, and conformal rescaling, simplifies  $L_{0}$ so that finally
\begin{eqnarray}
\tK &=& \p_{t}, \\
L_{0}&=& \p_{t} - \beta (x,y,z) \p_{z} + \lambda \p_{z}, \label{L0} \\
L_{1}&=&\partial_x + \lambda \p_{y} + (A_0(x,y) + \lam A_1(x,y) + \lam^2
A_2(x,y) + \lam^3 A_3(x,y))\p_{\lam} \label{L1} \\ && +
(C(x,y,z) + \lam D(x,y,z))\p_t + (E(x,y,z) + \lam F(x,y,z)) 
\p_z \nonumber.
\end{eqnarray}

One can read off a metric $g \in [g]$ corresponding to the twistor
distribution given by (\ref{L0}) and (\ref{L1}) by comparing with
(\ref{projtwistordistribution1}) and (\ref{projtwistordistribution2})
and reading off a null tetrad. One finds that $\bK \wedge d \bK =
\beta_z dx \wedge dy \wedge dz$, where $\bK = g(\p_t,.)$.
Thus the twist of the Killing vector $\p_t$ vanishes iff $\beta$ does
not depend on $z$. Since existence of twist is a conformally invariant
property, the cases $\beta_z = 0$ and $\beta_z \neq 0$ are
genuinely distinct, not an artefact of our coordinate choices. We now
analyse each in turn.

\vspace{10pt}
\textbf{Twist-free case:} $\beta_z = 0.$ 
Calculating the commutator of $L_0$ and $L_1$ we obtain
\begin{eqnarray} \label{nontwistcomm1}
[L_0,L_1]=&& (-\beta + \lambda)(C_z + \lambda D_z)\p_t + (\beta_x +
\lambda \beta_y - \beta E_z - \lambda \beta F_z + \lambda E_z +
\lambda^2 F_z - \nonumber\\  
&&(A_0+\lambda A_1 + \lambda^2 A_2 + \lambda^3 A_3))\p_z.
\end{eqnarray}
Since we require $\{ L_0, L_1 \}$ to be integrable, this must be a
multiple of $L_0$. We deduce that \footnote{In \cite{calderbank1} the
  resulting equations are interpreted as a special case of a gauge theory
  defined on a projective surface. A solution is called a
  \emph{projective Higgs pair}. This also applies to the twisting case.}
\be \label{nontwistcomm2} 
[L_0,L_1] = (-\beta + \lambda)(C_z + \lambda D_z) L_0.
\ee
Now comparing the $\p_z$ coefficients of (\ref{nontwistcomm1}) and
(\ref{nontwistcomm2}) we get four equations, one for each power of
$\lambda$. We can solve three of them, and use $L_1\rightarrow L_1-CL_0$
which does not change the
conformal structure. This yields

%The $\lambda^3$ equation is
%\be \label{Deqn}
%D_z = -A_3 \rightarrow D(x,y,z) = - z A_3(x,y) + Q(x,y),
%\ee
%where $Q$ is arbitrary. 

%The $\lambda^2$ equation is 
%\be \label{notwist2}
%F_z - C_z = A_2 + 2 \beta A_3.
%\ee
%Substituting (\ref{Deqn}) into this and solving gives
%\be \label{Feqn}
%F(x,y,z) - C(x,y,z) = z(A_2(x,y) + 2 \beta(x,y)A_3(x,y)) + P(x,y),
%\ee
%where $P$ is arbitrary.

%The $\lambda$ equation is
%\be \label{notwist1}
%E_z + \beta C_z = - \beta_y + A_1 + \beta A_2 + \beta^2 A_3,  
%\ee
%which integrates to gives
%\be \label{Eeqn}
%E(x,y,z) + \beta(x,y) C(x,y,z) = z (-\beta_y + A_1 + \beta A_2 +
% \beta^2 A_3) + R(x,y),
%\ee
%where $R$ is arbitrary. Substituting $D,F,E$ from (\ref{Deqn}),
%(\ref{Feqn}), (\ref{Eeqn})
%into $L_1$, we find the following:
\begin{multline} \label{L1minusCL0}
L_1 = \partial_x + \lambda \p_{y} + (A_0 + \lam A_1 + \lam^2
A_2 + \lam^3 A_3)\p_{\lam} \\ +
\lambda(-z A_3 + Q)\p_t + (z(-\beta_y + A_1 + \beta A_2 + \beta^2 A_3) +
+ \lambda(z(A_2 + 2 \beta A_3) + P)) \p_z, 
\end{multline}
where $P$ and $Q$ are arbitrary functions of $(x, y)$ and we have
eliminated one arbitrary function by translating the $z$ coordinate.
There is
one remaining equation to solve, corresponding to the
$\lambda^0$ coefficient of $\p_z$. This equation is as follows:
\be \label{Keqn}
\beta_x + \beta \beta_y - A_0 - \beta A_1 - \beta^2 A_2 - \beta^3 A_3 = 0.
\ee
The metric (\ref{nontwistinggeneral}) in Theorem \ref{localform}
corresponds to the twistor distribution given by $L_0$, with $\beta_z =
0$, and (\ref{L1minusCL0}).
If $\beta(x,y)$ is regarded as defining a section of
$PTU$, then (\ref{Keqn}) says that this section is tangent to lifted
geodesics of the projective structure. In terms of the base, a
solution is given by a congruence of geodesics.

\vspace{10pt}
\textbf{Twisting case:} $\beta_z \neq 0.$ We may perform a coordinate transformation $z
\rightarrow  \beta
(x,y,z)$. This does not affect the general form (\ref{L1}) of
$L_1$. Performing the coordinate change and dividing by $\beta_z$
gives the following form for $L_0$:
\be \label{L0withomega}
L_0 = H(x,y,z) \p_t - z \p_z + \lam \p_z,
\ee
where $H$ is a non-zero arbitrary function. Calculating the 
commutator gives
\begin{eqnarray*}
[L_{0},L_{1}] &=& ((-z+\lambda)(C_{z}+\lambda D_{z})-(E+\lambda 
F)H_{z})\p_{t}+ \\ && ((-z+\lambda)(E_{z}+\lambda F_{z})-(E+\lambda 
F)-(A_{0}+\lambda A_{1} + \lambda^{2} A_2 + \lambda^{3} A_{3})) 
\p_{z}.
\end{eqnarray*}    
We require $[L_{0},L_{1}]=\alpha L_{0}$ for some function $\alpha(x,y,z,\lambda)$, which is at most quadratic in $\lambda$, 
since otherwise $\alpha L_{0}$ will contain powers of $\lambda$ 
greater than three, and such terms do not occur in the commutator
above. 
We  make a replacement 
$L_1\rightarrow L_1-FL_0$, and analyze equations obtained from
comparing 
coefficients of  $\p_{z}, \p_{t}$. 
This puts $L_1$ in the form
\begin{eqnarray*} 
\label{L1intermediate}
L_{1}&=&\partial_x + \lambda \p_{y} + (A_0 + \lam A_1 + \lam^2
A_2 + \lam^3 A_3)\p_{\lam} \\ && +
(C + \lam D)\p_t + (A_{0} + z A_{1} + z^{2} A_{2} + 
z^{3} A_{3}) \p_{z},
\end{eqnarray*}
where $C(x,y,z), D(x,y,z), H(x,y,z)$ satisfy
\begin{eqnarray}
C_{z}-2zD_{z} &=& -H A_{2} + H_{y}. \label{lambda2}\\
D_{z} &=& -H A_{3}. 
\end{eqnarray}
and
\be \label{omegasprayeqn}
(\p_{x} + z \p_{y} + (A_{0} + z A_{1} + z^{2} A_{2} + 
z^{3} A_{3})\p_{z}) H = 0.
\ee
The only things remaining now are to find expressions for $C$ and $D$ 
and construct the metric. In order 
to integrate equations (\ref{lambda2}) it is 
convenient to express $H(x,y,z)$ as the second derivative of 
another function $G(x,y,z)$, i.e. we set
$$
H(x,y,z) = \frac{\p^{2} G}{\p z^{2}}(x,y,z).
$$
Then equations (\ref{lambda2}) integrate to give
\begin{eqnarray*}
C &=& -G_{z} A_{2} - 2 A_{3} (z G_{z} - G) + G_{zy} + 
\rho(x,y), \\
D &=& - G_{z} A_{3} + \sigma(x,y),
\end{eqnarray*}
where $\rho$ and $\sigma$ are arbitrary functions. 
Notice that $G$ has a `gauge freedom' $G \rightarrow G + z \gamma(x,y) +
\delta(x,y)$, since
(\ref{Gequation}) will still be satisfied. Using this and a coordinate
change $t \rightarrow t + \xi(x,y)$, one can set
the functions $\rho$ and $\sigma$ to zero.
The twistor 
distribution $\{ L_{0}, L_{1} \}$ is now fully determined:
\begin{eqnarray*}
L_0 &=& G_{zz} \p_t - z \p_z + \lambda \p_z, \\
L_1 &=& \p_x + \lambda \p_y + (A_0 + \lam A_1 + \lam^2 A_2 + \lam^3
A_3)\p_{\lam} \\ && + (-G_z A_2 - 2 A_3 (z G_z - G) + G_{zy} ) -
 \lam (G_z A_3)) \p_t \\ && + (A_{0} + z A_{1} +
z^{2} A_{2} + z^{3} A_{3}) \p_z.
\end{eqnarray*}
The distribution is integrable iff $G$ satisfies (\ref{Gequation}).
Calculating the corresponding null tetrad gives the conformal
structure (\ref{twistinggeneral}) in Theorem \ref{localform}.$\Box$

\section{Examples} \label{examples}

\subsection{Neutral Fefferman conformal metrics.} If $G_{zz}$ is simply a constant, then
(\ref{Gequation}) is satisfied. So given any projective structure and
setting $G_{zz}=1$ we obtain a family of conformal structures with
twist which
reduce to the given projective structure. Solving for $G$ gives
$$
G = \frac{z^2}{2} + z \gamma(x,y) + \delta(x,y).
$$
The corresponding metric takes the form
\begin{multline} \label{specialtwisting}
(dt + ((z+\gamma)A_3 + \sigma)dy + ((z+\gamma)A_2 + 2A_3(\frac{z^2}{2}
- \delta) - \gamma_y + \rho)dx)(dy - z dx) \\ - (dz - (A_0 + z A_1 +
z^2 A_2 + z^3 A_3) dx) dx,
\end{multline}
where we have chosen not to eliminate $\sigma$ and $\rho$.
By direct calculation one can show that the ASD Weyl tensor has
Petrov-Penrose type III or N, and it is type N precisely when the
following hold:
\begin{eqnarray*}
\gamma A_3 + \sigma &=& \frac{1}{3} A_2, \\
\gamma A_2 - 2 A_3 \delta - \gamma_y + \rho &=& \frac{2}{3} A_1.
\end{eqnarray*}
One can always choose $\rho, \sigma, \gamma, \delta$ so that these
are satisfied. In this case, the metric is the same as (31) in
\cite{nurowskisparling}, with their $Q$ cubic in $p$. These are
neutral signature analogues of Fefferman metrics.

\subsection{ASD pp-waves} Notice that the metric (\ref{nontwistinggeneral})
does not explicitly contain the function $A_0(x,y)$ of the projective
structure. The metric is always ASD for any choice of
$\beta,A_1,A_2,A_3$; one can regard (\ref{Keqn}) as
giving $A_0(x,y)$ in terms of these functions. One the other hand, if
one wants to specify $A_0$, then one must choose a solution of (\ref{Keqn}) for
$\beta$. In the special case $A_0=0$, we have the solution
$\beta=0$. One then obtains the following metric:
\be \label{betazero}
g = (dt + (P+z A_2) dx + (Q+z A_3 dy) ) dy - (dz + z A_1 dx ) dx.
\ee

Different choices of function $\beta(x,y)$ in
(\ref{nontwistinggeneral}) can give rise to different metrics. Suppose
we choose the flat projective structure. Then
$\beta(x,y)$ must satisfy the equation (\ref{beta_equation})
with $A_\alpha=0$.
By direct calculation one can show that the metric
(\ref{nontwistinggeneral}) is type III iff $\beta_{yy} \neq 0$,
otherwise it is type N. So the conformal structures  
with $\beta_{yy}=0$ and
$\beta_{yy}\neq 0$ are genuinely distinct.

\subsection{Pseudo-hyper-K\"ahler metrics} \label{PHK}
We will find some examples of neutral ASD metrics with
null conformal Killing vectors by independent means, and interpret them using
our results. We will use Pleba\'nski's method \cite{plebanski} adapted
to neutral signature, which converts the
problem of finding Ricci-flat ASD neutral metrics, or
\emph{pseudo-hyper-K\"ahler}, to the problem of
solving a non-linear second order PDE. He showed
that such metrics are locally of the form
\begin{equation} \label{hyperkahler}
g = dY (dT - \Theta_{XX} dY - \Theta_{TX} dZ) - dZ (dX + \Theta_{TT} dZ
+ \Theta_{TX} dY),
\end{equation}
where $\Theta(T,X,Y,Z)$ satisfies the `second Heavenly Equation':
\begin{equation} \label{heavenly}
\Theta_{YT} - \Theta_{ZX} + \Theta_{TT} \Theta_{XX} - \Theta_{XT}^2 = 0.
\end{equation}
The primed connection coefficients vanish when using the tetrad indicated in
(\ref{hyperkahler}), so there is a basis of covariantly constant
primed spinors $o^{A'}=(1,0)$, $\iota^{A'} = (0,-1)$. There is therefore also
a basis $\Sm^{A'B'}$ of covariantly constant null self-dual
two forms, written in spinors as follows:
\begin{eqnarray}
\Sm^{0'0'} &=& \frac{1}{2} \ \iota_{A'} \iota_{B'} \ep_{AB} \ \mathbf{\theta}^{AA'} \wedge
\mathbf{\theta}^{BB'},
\label{sig00}
\\
\Sm^{0'1'} = \Sm^{1'0'} &=& \frac{1}{2} \ o_{(A'} \iota_{B')} \ep_{AB} \mathbf{\theta}^{AA'}
\wedge \mathbf{\theta}^{BB'}, \label{sig01}
\\
\Sm^{1'1'} &=& \frac{1}{2} \ o_{A'} o_{B'} \ep_{AB} \ \ \mathbf{\theta}^{AA'} \wedge
\mathbf{\theta}^{BB'}. \label{sig11}
\end{eqnarray}
Using the identification between two-forms and endomorphisms given by
$g$, we can write
$$
R = \Sm^{0'0'} - \Sm^{1'1'}, \ \ \ I = \Sm^{0'0} + \Sm^{1'1}, \ \ \ S
= \Sm^{0'1'}.
$$
As endomorphisms, these satisfy
\begin{equation} \label{endos}
-I^2 = R^2 = S^2 = \text{Id}, \ \ IRS = \text{Id},
\end{equation}
which is easy to check using their spinor forms. There is a
hyperboloid's worth of almost complex structures, $aI + bR + cS$, where $a^2
- b^2 - c^2 = 1$, which are parallel and hence integrable. This is a
\emph{pseudo-hyper-K\"ahler} structure.

Now writing (\ref{killingidentity}) using spinors by means of
(\ref{killing}) and (\ref{riemann}) gives
\begin{equation*}
\iota^A o^{A'} C_{ABCD} \ep_{A'B'} \ep_{C'D'} \\ = \nabla_{BB'}
(\phi_{C'D'} \ep_{CD} + \psi_{CD} \ep_{C'D'} + \frac{1}{2} \eta
\ep_{CD} \ep_{C'D'}),
\end{equation*}
where we have used Ricci flatness and anti-self-duality. For a pure
Killing vector or a homothety ($\eta$ constant), it follows that
\begin{equation} \label{sdderivative}
\nabla_{AA'} \phi_{B'C'} = 0.
\end{equation}
Therefore $\phi_{B'C'}$
is actually constant in the basis shown in (\ref{hyperkahler}). Now
let us suppose we have a null Killing vector which preserves the
$\alpha$--plane distribution spanned by $o^{A'} \e_{AA'}$. Then $K = \iota^A
o^{A'}\e_{AA'}$, for some $\iota^A$ and using
(\ref{alg2}) and (\ref{sdderivative}) we get
$$
\phi_{B'C'} =  a_1 \ o_{B'} o_{C'}  + a_2 \ o_{(B'} \iota_{C')} ,
$$
for constant $a_1$, $a_2$. Consider the three distinct cases:
$\phi_{B'C'}$ vanishing ($a_1=a_2=0$), non-vanishing but degenerate ($a_1
\neq 0, a_2=0$), and non-degenerate ($a_1=0, a_2 \neq 0$). For
$K=\p_T$ we 
get the first case, $K=Y \p_X + Z \p_T$ the second,
and $T \p_T + X \p_X$ the third, and with some efford it can be shown
that  these choices are
canonical ( the first two cases were analysed in
\cite{finleyplebanski}). 
In order
for any of these to be
Killing, an equation for $\Theta$ coming from the Killing equation
must be satisfied. In fact we were only able to fully solve for the
first two cases. 

$\bullet$ \ $K=\partial_T$ \newline
Since $\p_T$ has no twist we expect this to be of the form
(\ref{nontwistinggeneral}). It is a neutral signature version of a
tri-holomorphic Killing vector; i.e. it Lie-derives $I,R,S$. Solving
the Killing equations in
conjunction with (\ref{heavenly}) results in the following metric:
\begin{equation} \label{case1}
g = dY dT - dZ dX - Q(X,Y) dY^2,
\end{equation}
where $Q$ is an arbitrary function. This is simply the
split-signature pp-wave metric, and is a special case of
(\ref{betazero}). Here $K$ is a self-dual Killing vector in the sense of
Gibbons et al. \cite{BGPPR}.

The local expression (\ref{case1}) in this example corresponds
to a class of global neutral metrics on compact four--manifolds. To see this
we compactify the flat projective space $\R^2$ to two--dimensional torus
$U=T^2$ with the projective structure coming from the flat metric. Both
$T$ and $Z$ in (\ref{case1})
 are taken to be periodic, thus leading
to $\hat{\pi}:M\longrightarrow U$, the holomorphic 
toric fibration over a torus. Assume the suitable 
periodicity on the function $Q:U\longrightarrow \R$. 
This leads to a commutative diagram
\begin{eqnarray*}
&M&\\
T^2&\downarrow&\searrow \hat{\pi}^*Q\\
&U& \stackrel{Q}{\longrightarrow}\R .
\end{eqnarray*}

This example can be put into the framework of 
\cite{kamada} and \cite{Fino1},
where $M$ is regarded as a primary Kodaira surface $\C^2/G$ and
$G$ is the fundamental group of $M$ represented injectively in
the group of complex affine transformations of $\C^2$.  In this framework
the K\"ahler structure on $M$ is given by 
$\Omega_{flat}+i\p\overline{\p}(\hat{\pi}^*Q)$, where $(\p, \Omega_{flat})$ 
is the flat K\"ahler structure on the Kodaira surface induced from $\C^2$.

$\bullet$ \ $K = Y \partial_X + Z \partial_T$ \newline
Again, this is twist-free and we expect the metric to be of the form
(\ref{nontwistinggeneral}). Solving the Killing equations in
conjunction with (\ref{heavenly}) results in the following metric:
\begin{equation} \label{case2}
g = dY dT - dZ dX
-\frac{H(\frac{Y}{YT-ZX},\frac{Z}{YT-ZX})}{(YT-ZX)^3}(Y dZ - Z dY)^2,
\end{equation}
where $H$ is an arbitrary analytic function of two variables. This is
a generalization of the Sparling-Tod metric \cite{sparlingtod}. It is
easy to show that the arguments of $H$ are in fact constant on the
special \bs s, so serve as coordinates on $U$. 

%Now
%perform the following coordinate transformation:
%\begin{eqnarray*}
%t &=& \frac{1}{2}(\frac{X}{Y} - \frac{T}{Z}), \\
%z &=& \frac{1}{2} \log(YZ),\\
%x &=& \frac{Z}{YT+ZX}, \\
%y &=& \frac{Y}{YT+ZX}.
%\end{eqnarray*}
%The main point of this transformation is that the Killing vector
%becomes $\p_t$. After a conformal transformation and a
%redefinition of the arbitrary function, the metric
%becomes:
Using the following coordinate transformation
\begin{eqnarray*}
t &=& -\frac{1}{2}(\frac{X}{Y}+\frac{T}{Z}),\\
z &=&  (YZ)^{-\frac{1}{2}},\\
x &=& \frac{YT-XZ}{(YZ)^{\frac{1}{2}}},\\
y &=& \text{log} \big( \frac{Z}{Y} \big).
\end{eqnarray*}
the metric (\ref{case2}) takes the following form:
$$
g = \frac{1}{z^2}(dy dt - dz dx + z
A_3(x,y) dy^2),
$$
where now the Killing vector is $\p_t$. Multiplying by the conformal
factor $z^2$, we get a special case of
(\ref{betazero}).
The projective structure is non-trivial, unlike for
the pp-wave above. The projective structure is
special in that it depends on only one arbitrary function.

$\bullet$ \  $ T \partial_T + X \partial_X$ \newline
In this case we were not able to fully solve the Killing equations in
conjunction with (\ref{heavenly}). This Killing vector is twisting, so
the answer must be of the form (\ref{twistinggeneral}).

%\section{Special classes of metric} \label{specialmetrics}
%Having obtained our local classification, it is interesting to look at
%the metrics in more detail. One can ask questions like: What is the
%Petrov classification of the Weyl curvature? Which conformal classes
%contain Ricci-flat metrics or, more generally, Einstein metrics? Are
%there conformal classes which do not contain any such metrics?

%Determining the Petrov classification is a straightforward algebraic
%question. The other questions are more subtle. There
%are general methods available to tackle whether a conformal class
%contains an Einstein metric \cite{szekeres}, but often the
%computations involve such long expressions as to make them
%impractical, even by computer. We have adopted a more `hands-on'
%approach, where one simply puts in a conformal factor and looks at the
%equations that give Ricci-flatness. 

\subsection{Pseudo-hyper-hermitian conformal structures}
This is a generalization of the pseudo-hyper-K\"ahler case discussed
in the last section. We will refer to a neutral metric $g$ as
pseudo-hyper-hermitian (also called hyper-para---hermitian \cite{ivanov})
when there exist endomorphisms $I,R,S$
satisfying the algebra (\ref{endos}), such that any complex structure
$\mathcal{J}_{(a,b,c)} = a I + b R + c S$ is integrable for $a^2 - b^2
-c^2 = 1$, and $g$ is hermitian with respect to any of these complex
structures. For $g$ to be hermitian with respect to a complex
structure $\mathcal{J}$ means $g(\mathcal{J} X, \mathcal{J} Y) =
g(X,Y)$. Note that for pseudo-hyper-K\"ahler, the endomorphisms
$I,R,S$ must also be covariantly constant with respect to the
Levi-Civita connection of $g$.

In \cite{dunajski1}, it is shown that one can always find a tetrad for
a pseudo-hyper-hermitian metric such that the twistor distribution has
no $\p_{\lambda}$ terms. Equivalently, the twistor space fibres over
$\CP$. Now let us suppose that we have a null conformal Killing that
is tri-holomorphic, i.e. it preserves $I$, $R$ and $S$ and so it
preserves the holomorphic fibration ${\cal PT}\rightarrow \CP$.
All such cases
are classified by the following
\begin{proposition}
All pseudo-hyper-hermitian ASD metrics with triholomorphic null conformal
Killing vectors are of the form (\ref{nontwistinggeneral}) or
(\ref{twistinggeneral}) up to a conformal factor, where the corresponding ODE
(\ref{secondorderode}) is point equivalent to a derivative of a first
order ODE.
\end{proposition}
\emph{Proof.} Let $g$ be a pseudo-hyper-hermitian ASD metric, and $K$
be a triholomorphic conformal Killing vector. Since $g$ is ASD, it
follows from Theorem \ref{localform} that there are coordinates such
that, up to a conformal factor, $g$ is of the form
(\ref{nontwistinggeneral}) or (\ref{twistinggeneral}). From
\cite{dunajski1}, it is possible to find a tetrad such that the
twistor distribution has no $\p_{\lambda}$ terms. Now a change in
tetrad corresponds to a M\"obius transformation of
$\lambda$. Since $K$ is triholomorphic, its lift will have no
$\p_{\lambda}$ terms in the tetrad where the twistor distribution has
no $\p_{\lambda}$ terms. Therefore the M\"obius
transform does not depend on $t$, otherwise $\p_t$ will no longer
Lie-derive the twistor distribution (one would have to add
$\p_{\lambda}$ terms). Furthermore, the M\"obius
transformation does not depend on $z$, otherwise $\p_{\lambda}$ terms
will be introduced into $L_0$. Hence there is a M\"obius
transformation of $\lambda$, depending only on $(x,y)$, such that the
$\p_{\lambda}$ terms in $L_1$ are eliminated.

After this change in $\lambda$, the projective structure spray in
$L_1$ will be of the following form:
$$
\cS = a \p_x + b \p_y + \lambda (c \p_x + e \p_y),
$$
where $a,b,c,e$ are functions of $(x,y)$ with $ae-bc \neq
0$. Coordinate freedom
$(x,y) \rightarrow (\hat{x}(x,y), \hat{y}(x,y))$ and scaling freedom
(the projective structure is unchanged if $\cS$ is multiplied by a
non-zero function) allows us to set $a=1$, $c=0$, $e=1$, giving
$\cS = \p_x + (b + \lambda) \p_y$. Now perform another M\"obius
transformation $\lambda \rightarrow b + \lambda$, which gives the
following spray:
\begin{equation} \label{specialprojspray}
\p_x + \lambda \p_y + (b_x + \lambda b_y) \p_{\lambda}.
\end{equation}
This corresponds to the second-order ODE 
\begin{equation} \label{specialODE}
\frac{d^2 y}{dx^2} = A_1 (x,y) \Big( \frac{dy}{dx}
\Big) + A_0(x,y),
\end{equation}
where $A_1 = \frac{\p b}{\p y}$, $A_0 = \frac{\p b}{\p x}$ for a
function $b(x,y)$.
This is the derivative of the general first-order ODE
\begin{equation}
\frac{dy}{dx} = b(x,y).
\end{equation}
Hence the original projective structure is point-equivalent to the one
corresponding to (\ref{specialODE}).
$\square$
\vspace{10pt}

Note that if a (holomorphic) projective structure spray contains no
$\p_{\lambda}$ terms, its twistor space fibres over $\CP$, since each
integral curve can be labelled by the $\lambda$ coordinate. So a
by-product of the proof of the above Proposition and Theorem
\ref{hitchin_lebrun}
is the following
\begin{proposition} \label{specialprojective}
There is a one to one correspondence between holomorphic
2D projective structures s.t. the corresponding second order
ODE is point equivalent to the
derivative of a first order ODE, 
and complex surfaces which contain a holomorphic curve with normal
bundle ${\cal O}(1)$ and fiber holomorphically over $\CP$.
\end{proposition}
This is of interest purely as a statement about
projective structures. Note that although all first order ODEs can be
transformed to the trivial first order ODE $d y / d x = 0$ by
coordinate transformation, this does not mean that the derivative of
any such equation is flat, in the sense of Section
\ref{flatness}. This can be shown by calculating the invariant
(\ref{flatprojstructure}) for (\ref{specialODE}) and showing that it
does not necessarily vanish.

\subsection{Conformal structures containing no Ricci-flat metrics}
\label{confvac}
In this section we show that there are conformal structures of the form
(\ref{nontwistinggeneral}) which do not contain Ricci-flat
metrics. Before doing so we discuss the Petrov-Penrose classification
for the conformal structures (\ref{nontwistinggeneral}) and
(\ref{twistinggeneral}).
\begin{proposition}
Let $K^{AA'}=\iota^A o^{A'}$ be a null conformal Killing vector for
ASD conformal structure. Then $\iota^A$ is  a principal direction,
that is 
\be
\label{principal_dir}
\iota^A \iota^B \iota^C \iota^D C_{ABCD} = 0.
\ee
Moreover if the twist of $K$ vanishes 
the conformal structure  is of type $III$ or $N$, that is
\be
\label{III_N}
\iota^A\iota^B C_{ABCD}=0.
\ee
\end{proposition}

\emph{Proof}. From (\ref{alg1}) we have 
$$
\nabla_{AA'}(\iota^C \iota^D \psi_{CD}) = 0.
$$ 
Expanding this out we obtain
\be \label{inabpsi}
\iota^B \iota^C \nabla_{A A'}\psi_{BC} = -2 \psi_{BC} \iota^C
\nabla_{AA'} \iota^B = \iota_A \mu_{A'},
\ee
for some spinor $\mu_{A'}$. The last equality follows from
$(\ref{alg1})$ and $(\ref{GSF1})$.

Now pick a conformal frame in which $K$ is a pure Killing vector. The
well known identity $\nabla_a \nabla_b K_c = R_{bcad} K^d$ implies
$$
\nabla^{A'}_{\ A} \psi_{BC} = -2 C^D_{\ ABC}K^{A'}_D - 2
K^{B'}_{(A}\Phi^{A'}_{\ BC)B'}+\frac{1}{6}R \epsilon_{A(B}K^{A'}_{\ C)}
- \frac{4}{3}\epsilon_{A(B}\Phi_{C)}^{\ D D' A'}K_{DD'}.
$$
On contracting both sides by $\iota^A \iota^B \iota^C$ and using
(\ref{inabpsi}), all terms vanish except the term involving $C^D_{\
  ABC}$, giving (\ref{principal_dir}).

Now let us assume that $K$ is non--twisting, i. e.  $\bK
\wedge d \bK = 0$ where $\bK := g(K,)$. The Frobenius theorem implies
the existence of functions $P$ and $Q$ such that $\bK =PdQ$. We can
now choose a conformal
factor such that $d \bK = 0$. Then $K$ is covariantly constant
($\nabla_a K_b = 0$), and we deduce
\begin{eqnarray}
\nabla_{AA'} \iota_B &=& A_{AA'} \iota_B, \label{nabiota} \\
\nabla_{AA'} o_{B'} &=& -A_{AA'} o_{B'}, \label{nabmicron}
\end{eqnarray}
for some one-form $A_{AA'}$. Consider the spinor Ricci identity
\cite{penroserindler}
$$
\bigtriangleup_{A'B'}o_{C'} = (C_{A'B'C'D'}-\frac{1}{12} R
\epsilon_{D'(A'}\epsilon_{B')C'})o^{D'},
$$
where $\bigtriangleup_{A'B'} = \nabla_{A
  (A'}\nabla_{B')}^A$. Substituting ($\ref{nabmicron}$) into this and
using $C_{A'B'C'D'}=0$ gives
$$
o_{C'}\nabla_{A(A'}A^A_{B')} = -\frac{1}{12} R
o_{(A'}\epsilon_{B')C'}.
$$
By contracting with $o^{C'}$ we find $R=0$. Now consider the Ricci identity
$$
\bigtriangleup_{AB}\iota_{C} = (C_{ABCD}-\frac{1}{12} R
\epsilon_{D(A}\epsilon_{B)C})\iota^{D}.
$$
Substituting $R=0$ and (\ref{nabiota}) into this gives
$$
\iota_C \nabla_{A'(A}A^{A'}_{B)} = C_{ABCD}\iota^D.
$$
Contracting this with $\iota^C$ gives (\ref{III_N}),
from which it follows that the curvature is type III or N. $\square$
\vspace{12pt}

In the twisting case the algebraic type of the Weyl spinor can be
general.  This can be shown by using the following two scalar invariants
\cite{penroserindler}:
$$
I = C_{ABCD} C^{ABCD}, \ \ \ \ J = C_{AB}^{\ \ CD}C_{CD}^{\ \ EF}
C_{EF}^{\ \ AB}.
$$
The condition for type III is $I=J=0$, and for type II that $I^3 =
6J^2$. Now consider the metric (\ref{twistinggeneral}), with the flat
projective structure $A_i=0, i=0,\ldots,3$. The function $G_{zz}$
satisfies
$$
(\p_x + z \p_y)G_{zz} = 0,
$$
which is solved in general when $G_{zz}$ is an arbitrary function of
$(zx-y)$. Suppose $G$ is given by:
$$
G(x,y,z) = \frac{\e^{zx-y}}{x^2} + z B(x,y),
$$
where $B(x,y)$ is arbitrary, so $G_{zz} = \e^{zx-y}$. Then the two scalar
invariants are as follows:
\begin{eqnarray}
I &=& -\frac{3}{2} x B_{yy} \e^{-3(zx-y)}, \\
J &=& \frac{3}{8}x(xB_{yyx}+3B_{yy}+xzB_{yyy}) \e^{-4(zx-y)}.
\end{eqnarray}
Therefore, from the conditions above, the metric is neither type II
nor type III. 

%Define the following two scalar invariants:
%$$
%I = C_{ABCD} C^{ABCD}, \ \ \ J = C_{AB}^{\ \ \ CD}C_{CD}^{\ \ \
%  EF}C_{EF}^{\ \ \ AB}.
%$$
%It is shown in \cite{penroserindler} that $I = J = 0$ implies type III
%or N. By direct calculation one finds these are satisfied for
%(\ref{nontwistinggeneral}), hence any ASD metric with null nontwisting
%conformal Killing vector is type III or N. 

To find metrics that are not conformally Ricci-flat we
use results of Szekeres \cite{szekeres}. Although these were
derived for Lorentzian signature, they can also be applied to our ASD
neutral signature case, essentially because the Weyl curvature is
still made up of a single spinor $C_{abcd}=C_{ABCD} \ep_{A'B'}
\ep_{C'D'}$ as in the Lorentzian case (of course in Lorentzian case it
is complex hermitian, not real).

Consider the metric
(\ref{betazero}) with $A_1=0$.
By direct calculation, one finds that
$C_{ABCD}$ is type N iff $(A_2)_x = 0$, otherwise it is type
III.  Now suppose (\ref{nontwistinggeneral}) is type III, i.e. 
$(A_2)_x \neq 0$. The reason for this
is that we can apply a result of Szekeres to obtain an obstruction to
Ricci-flatness. It is shown in $\cite{szekeres}$ that for types I, II, D or
III, a necessary condition for existence of a Ricci-flat metric in the
conformal class is the following tensor equation
$$
-\frac{1}{2}C_{pqfh} C_{rs}^{\ \ f h} C_{a b c \ ;d}^{\ \ \ d} +
(C_{pq}^{\ \ df} C_{rsf \  ;h}^{\ \ \ h}
+C_{rs}^{\ \ df} C_{pqf \  ;h}^{\ \ \ h} ) = 0.
$$
This is just the tensor version of the spinor identity (3.1), page 209
\cite{szekeres}. Calculating this one finds that
$(A_2)_{xx}$ is an obstruction to its vanishing (we used MAPLE for the
calculation). Therefore we have a
class of non-conformally vacuum type III neutral ASD conformal structures with
non-twisting null conformal Killing vectors. 
%%%%%%%%%%%%%%%%%%%%%%%%%%%%%%%%%%%%%%%%%%%%%%%%%%%%%%%%%%%%%%%%%%%%%
\section{Twistor reconstruction} \label{twistorreconstruction}
\setcounter{equation}{0}
We have shown that when a conformal structure $[g]$ has a null
conformal Killing vector, the twistor space $\PT$ fibres over the
twistor space of a projective structure, and we have classified the
possible local forms for such conformal structures.

The twistor lines in a projective structure twistor space $Z$ have
normal bundle $\OO(1)$. The twistor lines in a conformal structure
twistor space have normal bundle $\OO(1) \oplus \OO(1)$. Let
$\mathcal{B}$ be a holomorphic fibre bundle over $Z$ with one
dimensional fibres. Let $\hat{u}$ be a twistor line in $Z$. Then if we
want $\mathcal{B}$ to be a conformal structure twistor space, the
normal bundle of $\hat{u}$ in $\mathcal{B} |_{\hat{u}}$ must be $\OO(1)$.
Given a projective structure twistor space, one  way of forming
a fibre bundle with the correct property is to take a power of the
canonical bundle $\kappa$, which reduces to $\OO(-3)$ on twistor
lines. The bundle $\kappa^{-1/3}$ reduces to $\OO(1)$ on twistor
lines, and exists provided we take $Z$ to be a suitably small
neighbourhood of a twistor line. So the total space of $\kappa^{-1/3}$
is a conformal structure twistor space. 

Consider the simplest possible case, where $Z$ is the total space of
$\OO(1)$, corresponding to a flat projective structure. In this case
$\kappa^{-1/3}$ is the total space of $\OO(1) \oplus \OO(1)$, the
twistor space of the flat conformal structure.
To go further, note that given a line bundle $\mathcal{A}$ over $Z$
which reduces to $\OO(1)$ on twistor lines, any affine bundle modelled
on $\mathcal{A}$ will also have the correct property on twistor
lines. In the simplest case described above, taking affine bundles
modelled on $\kappa^{-1/3}$ results in the the twistor space of the
pp-wave metric (\ref{case1}). In fact, this is precisely the first case
discussed by Ward in $\cite{ward1}$, although he does not phrase it in
this way. We will now show how this works.

\subsection{Example 1. PP-waves.}

First we will give a twistorial demonstration of a fact shown
in Section \ref{PHK}, namely that for a pseudo-hyper-K\"ahler metric with
triholomorphic null Killing vector $K=\iota^A o^{A'} \e_{AA'}$ with $o^{A'}$
covariantly constant, the resulting projective structure is flat. The
twistor space of an analytic pseudo-hyper-K\"ahler metric fibres over $\CP$,
$\sigma: \PT \rightarrow \CP$ \cite{penrose,hitchin}. There is a
section $\varpi$ of $\Lambda^2 \PT \times \sigma^* \OO(2)$. This is a
symplectic form of `degree 2' on the fibres. In the spin bundle
picture, $\varpi$ is the push forward to $\PT$ of the symplectic form
$\Sm = \Sm^{A'B'} \pi_{A'} \pi_{B'}$ on $S'$, where $\Sm^{A'B'}$ are defined
as in Section \ref{PHK}. This form is Lie-derived over the twistor
distribution as a consequence of the $\Sm^{A'B'}$ being covariantly
constant, and is homogeneous in the $\pi_{A'}$, so the push-forward is
well defined.

As explained in Section \ref{main}, $\cK$ vanishes on a hypersurface
$\mathcal{H}$ in $\PT$, where $\mathcal{H}$ is the projection to $\PT$
of the hypersurface $\pi.o = 0$ in $S'$. For $o^{A'}$ covariantly
constant, the function $1/(\pi.o)$ on $S'$ gives a section $\zeta$ of
$\sigma^* \OO(-1)$ on $\PT$, which blows up on $\mathcal{H}$. Then
$\zeta \otimes \cK$ is a non-vanishing `$\sigma^* \OO(1)$-valued' vector
field. Now in a local trivialization, $\zeta \otimes \cK$ Lie derives
the symplectic form $\varpi$, so it is Hamiltonian,
$$
\zeta \otimes \cK = \frac{\p h}{\p \omega^A} \frac{\p}{\p \omega_A},
$$
where $\varpi = \omega^0 \wedge \omega^1$. Now the $\omega^A$ should be
regarded as coordinates of `degree 1', that is they are coordinate
functions multiplied by a section of $\sigma^* \OO(1)$. Therefore for
the weights to agree, $h$ must be a section of $\kappa^* \OO(1)$,
rather than a bona fide
function. This gives a projection $\PT \rightarrow Z = \OO(1)$, with
fibres the trajectories of $\zeta \otimes \cK$, so the projective
structure twistor space is the total space of $\OO(1)$, which
corresponds to the flat projective structure.

Now suppose we start with the total space of $\OO(1)$ as the
minitwistor space $Z$. The twistor lines are
global holomorphic sections of
$\mathcal{O}(1)\rightarrow\CP$.

We will use a homogeneous coordinate description of $Z = \OO(1)$. Let
$\pi_{A'}$ be homogeneous coordinates for the base $\CP$
of $Z=\mathcal{O}(1)$, and let $\omega$ be a homogeneous coordinate for
the fibre of $Z=\OO(1)$. That is, $\OO(1) = \{ [\pi_{0'},\pi_{1'},\omega] :
[c\pi_{0'}, c \pi_{1'}, c\omega], c \in \mathbb{C}^*, [\pi_{0'},
  \pi_{1'}] \neq [0,0]   \}$.

Now cover the base $\CP$ in $\PT$ with two open sets $(\mathcal{U}_0, 
\mathcal{U}_1)$,
and lift this covering to $\mathcal{PT}$. Use homogeneous coordinates
$(\pi_{A'},\omega,\zeta_i)$ on $\mathcal{U}_i$.

The flat twistor space $\OO(1) \oplus \OO(1)$ can be formed as
follows. Consider the projection $\tau: \OO(1) \rightarrow \CP$. Then
$\OO(1) \oplus \OO(1)$ is the pull-back bundle $\tau^* \OO(1)$ over the
total space of $\OO(1)$. It is easy to check that this is the same as
taking $\kappa^{-1/3}$ where $\kappa$ is the canonical bundle of
$Z=\OO(1)$. To obtain curved twistor spaces, we can take
affine bundles over $\OO(1)$ modelled on $\tau^* \OO(1)$. To form these
we use the following transition functions:
\begin{displaymath}
\zeta_0=\zeta_1+f(\pi_{A'},\omega),
\end{displaymath}
where $f \in [f] \in H^1(Z,\tau^*\mathcal{O}(1))$, where $Z$ is
$\OO(1)$. The cohomology elements $f$ classify affine bundles over $Z$
modelled on $\tau^* \OO(1)$.

Global holomorphic sections of $Z \rightarrow
\CP$ are defined by $\omega = P(\pi_{A'})=\pi_{A'}x^{A'}$, with $x^{A'}=(X,Y)$
say. 

The sections of $\mathcal{PT}\rightarrow\CP$ are constructed
by putting $\zeta_i=\pi_{A'}t^{A'} + f_i$, where $t^{A'}=-(T,Z)$ say,
and $f=f_0-f_1$. The reason $f$ can be split in this way is
that when restricted to a twistor line in $Z$, $f$ becomes an element
of $H^1(\CP,\OO(1))$, and this group vanishes. To realise a splitting of $f$
we divide it by
$(\pi_{A'}o^{A'})^2$ for some constant $o^{A'}$, to get an element of
$H^1(Z,\tau^*\mathcal{O}(-1))$. Then we can use the fact that
$H^0(\CP,\OO(-1))=H^1(\CP, \mathcal{O}(-1))=0$, so any element can be
written as a difference of coboundaries, and the splitting is unique.
These sections are the $\CP$ twistor lines in $\PT$; we will refer to these
as $\hat{x}$, where $x$ is the point in $M$ with coordinates $(t^{A'}$, $x^{A'})$.

Let $\rho_{A'}$ be homogeneous coordinates on $\CP$.
The splitting is given by the Sparling formula:
\begin{displaymath}
\frac{f(\pi,P)}{(\pi.o)^2}=\oint_{\Gamma_0}\frac{f(\rho,P)}{(\rho.o)^2
\pi.\rho}\rho.d\rho-\oint_{\Gamma_1}\frac{f(\rho,P)}{(\rho.o)^2
\pi.\rho}\rho d\rho,
\end{displaymath}
where we are using Cauchy's integral formula, and $\Gamma_i \subset
\hat{x} \cong \CP$ are contours that bound a region
containing the point $\rho_{A'}=\pi_{A'}$. The
measure $\rho.d\rho$ means $\epsilon_{A'B'}\rho^{A'}d\rho^{B'}$.

Therefore
\begin{displaymath}
f_i=\oint_{\Gamma_i}\frac{(\pi.o)^2}{(\rho.o)^2}\frac{f(\rho,P)}{\pi.\rho}
\rho.d\rho. 
\end{displaymath}
The symplectic form $\varpi$ discussed above is given by
by $\varpi = d\omega  \wedge
d\zeta_i$ on $\mathcal{U}_i$. Restricting $\varpi$ to a section and taking
exterior derivatives keeping $\pi_{A'}$ constant, we obtain a formula
for $\Sm$, the pull-back of $\varpi$ to $S'$:
\begin{eqnarray*}
\Sm & = & d(\pi_{A'}x^{A'}) \wedge d(\pi_{B'}t^{B'}+f_0) \\
       & = & \pi_{A'} \pi_{B'} dx^{A'}\wedge dt^{B'} +
       \pi_{A'}dx^{A'}\wedge d f_0,
\end{eqnarray*}
where we are working over $\mathcal{U}_0$. Now

\begin{eqnarray*}
d f_0 &=& dx^{B'} \otimes \frac{\partial}{\partial x^{B'}} \oint_{\Gamma_0}
\frac{(\pi.o)^2}{(\rho.o)^2}\frac{f(\rho,\rho_{A'}x^{A'})}{\pi.\rho}\rho.d\rho\\
     &=& dx^{B'} \oint_{\Gamma_0} \frac{\rho_B'
  (\pi.o)^2}{(\rho.o)^2(\pi.\rho)}\frac{\partial f}{\partial P} \rho.d\rho.
 \end{eqnarray*}
Where we have used
$\frac{\partial}{\partial x^{A'}} \rightarrow \rho_{A'} \frac{\partial}{\partial
P}$. Using this we get

\begin{eqnarray*}
\Sm &=& \pi_{A'}\pi_{B'} dx^{A'} \wedge dt^{B'} +
\Big( \oint_{\Gamma_0}\frac{\pi_{A'}\rho_{B'} (o.\pi)^2}{(o.\pi)^2
  (\pi.\rho)}\frac{\partial f}{\partial P}\rho.d\rho\Big)
dx^{A'}\wedge dx^{B'} \\ 
       &=& \pi_{A'}\pi_{B'} dx^{A'} \wedge dt^{B'} + \frac{1}{2}\Big(
\oint_{\Gamma_0} \frac{(o.\pi)^2}{(o.\rho)^2}\frac{\partial
  f}{\partial P} \rho.d\rho \Big) dY \wedge dX \\
       &=& \pi_{A'}\pi_{B'} dx^{A'} \wedge dt^{B'} + (o.\pi)^2 Q(X,Y)
dY \wedge dX,
\end{eqnarray*}
where
\begin{displaymath}
Q(X,Y)=\frac{1}{2} \oint_{\Gamma_0}\frac{1}{(o.\rho)^2}\frac{\partial
  f}{\partial P} \rho.d\rho.
\end{displaymath}
Putting $o^{A'}=(1,0)$, we get the following
formula for $\Sm$ pulled back to $M\times\mathbb{C}^2$:
\begin{equation} \label{sigma}
\Sm= \pi_{0'}^2 (dT \wedge dX + Q(X,Y) dY \wedge dX) + \pi_{0'}
\pi_{1'} (dT \wedge dY - dX \wedge dZ) + \pi_{1'}^2 dZ \wedge dY.
\end{equation}
Calculating $\Sm$ in the Plebanski formalism from (\ref{sig00}),
(\ref{sig01}) and (\ref{sig11}) gives
\begin{eqnarray*}
\Sm &=& \pi_{0'}^2 (dT - \Theta_{XX} dY - \Theta_{TX} dZ)\wedge (dX +
\Theta_{TT} dZ + \Theta_{TX} dY)+ \\ &&\pi_{0'} \pi_{1'} (dT \wedge dY - dX
\wedge dZ) + \pi_{1'}^2 dZ \wedge dY.
\end{eqnarray*}
Comparing gives the forms $\Sm^{A'B'}$ and hence the  metric 
%\begin{displaymath}
%g=dYdT-dZdX-Q(X,Y)dY^2,
%\end{displaymath}
(\ref{case1}). The arbitrary function $Q$ corresponds to some
arbitrary cohomology element $f$.

\subsection{Example 2. Flat conformal structure.}
Here we show that given a conformal Killing vector for the flat
conformal structure, the underlying projective structure is also
flat. By the results of \cite{tafel}, we need only consider the conformal
Killing vectors $\p_T$ (non-twisting) and $T \p_T + X \p_X$
(twisting), where the flat metric is 
$$
g = dT dY - dX dZ.
$$
The non-twisting case is covered by the example of the last section,
with $Q(X,Y)=0$, so we know the projective structure is flat and
$Z=\OO(1)$. 

\begin{figure}
\caption{Quotient of the non--projective twistor space by the Euler vector 
field showing the singular set of ${\cal K}$.}
\label{case3_fig}
\begin{center}
\includegraphics[width=8cm,height=5cm,angle=0]{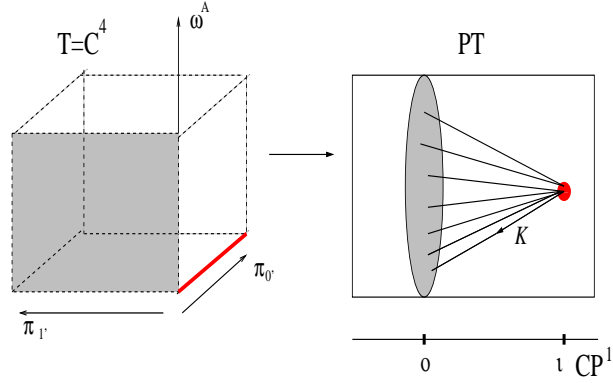}
\end{center}
\end{figure}

The twisting case is slightly more complicated. One can use the
spray picture, but instead we will analyse the twistor space $\PT$ and
show that the space of trajectories of $\cK$ is the flat projective
structure twistor space ${\mathbb{CP}^2}$
We work on  the non--projective twistor space ${\cal T}=\C^4$ with coordinates
$(\omega^A, \pi_{A'})$. The projective twistor space ${\cal PT}$ is a 
quotient of ${\cal T}$ be the Euler homogeneity vector field 
$\Upsilon=\omega^A/\p \omega^A+\pi_{A'}/\p \pi_{A'}$. The flat 
conformal class on the
complexified $\R^{2, 2}$ and the conformal twisting Killing vector are
represented by 
\[
g=\varepsilon_{AB}d p^Bd q^A, \qquad K=p^A\frac{\p}{\p p^A}
\]
where $x^{AA'}:=p^Ao^{A'}+q^A\iota^{A'}$ are coordinates on $M$. 
The point $(p^A, q^A)$ corresponds
to a two--plane in $\cal T$ given by solutions to the twistor equation
$\omega^A=x^{AA'}\pi_{A'}$.
The lift of $K$ to $S'$ is
\[
\tilde{K}=K+\pi_{1'}\frac{\p}{\p \pi_{1'}},
\]
and the orbits of the induced group action on the non--projective twistor space
are 
\[
\omega^A\rightarrow c\omega^A, \qquad \pi_{1'}\longrightarrow c\pi_{1'},\qquad
\pi_{0'}\longrightarrow \pi_{0'}.
\]
The holomorphic vector field on ${\cal T}$
\[
{\cal{K}}=\omega^A\frac{\p}{\p \omega^A}+\pi_{1'}\frac{\p}{\p \pi_{1'}},
\]
vanishes on the projective twistor space when it is proportional
to the Euler vector field. This happens on a set 
$B=\{\{\omega^A=0, \pi_{1'}=0\}\cup\{\pi_{0'}=0\}\}\subset{\cal T}$ 
which is a union of the line and
a hyperplane $\C^3\subset{\cal T}$.
The set $B$ descents to a union of a hypersurface and a point in the
projective twistor space (Fig. \ref{case3_fig}).
The minitwistor space $Z$ corresponding to the projective
structure $U$ is the factor space of ${\cal PT}/B$ by the trajectories of
${\cal K}$. Each trajectory in ${\cal T}$ is parametrised by its intersection
with the singular surface $\C^3$ given by $\pi\cdot o=0$ in ${\cal T}$
so the space of trajectories in ${\cal PT}$ is ${Z}={\mathbb{CP}^2}$. 
Two $\CP$s in ${\mathbb{CP}^2}$ intersect in
a point so the normal bundle of each $\CP$ is ${\cal O}(1)$ and we have
a projective structure. To obtain the explicit parametrisation of these 
$\CP$s eliminate $\pi_{0'}$ from the twistor equation 
to get $\pi_{1'}=\omega^Au_A$ where $u_A:=p_A/(p_Bq^B)$ parametrise the twistor
lines in ${Z}$ and are coordinates on $U$. The flat metric in 
$M$ is conformal  to  (\ref{twistinggeneral}) with $A_\alpha, G=z^2/2$ 
and conformal factor
$e^t$.

\section{Outlook} \label{conclusion}
We have locally classified neutral signature ASD conformal structures
with null Killing vectors. Some of these are defined on compact manifolds.
It would be interesting to  investigate the global properties of
other conformal structures we have found. 

It would also be interesting to understand in more detail which conformal
structures admit special types of metric, for example Ricci-flat or
Einstein (in this case the pure Killing vectors must be twist--free 
\cite{julia}). So far the only results we have in this direction are given
in Section \ref{confvac}. The existence of these  special 
metrics should be related to invariants of
the corresponding projective structure.

The recent work of Calderbank \cite{calderbank1} extended many of the
results obtained in this paper. In particular Calderbank gave a
twistor characterisation of ASD conformal structures which admit a
geodesic shear free congruence $\iota^A$. Not all such congruences
give rise to null conformal Killing vectors $K$ such that
$\iota^AK_{AA'}=0$, and Calderbank characterised those which do.

\vspace{10pt}

\textbf{Acknowledgements.} We wish to thank Helga Baum, 
David Calderbank, Claude LeBrun, Lionel Mason, George Sparling and Paul Tod 
for helpful discussions. S.W. thanks the EPSRC
for financial support.

\appendix

\section{Appendix} \label{appendix}
\setcounter{equation}{0}

Here we summarise the required spinor notation and present the
calculations leading to a proof of (\ref{Klift}).
We use similar conventions to Penrose and Rindler
\cite{penroserindler} adapted to neutral signature, but  our
indices are concrete.
\vskip5pt
{\bf Spin connection and curvature decomposition.}
As usual, we denote the Levi-Civita connection of the metric by
$\nabla$. The `spin connection coefficients' are defined by
$$
\nabla (\e_{CC'}) = \mathbf{\theta}^{DD'} \ot (\Gamma_{DD'C}^{\ \ \ \ \ \
  E}\e_{EC'} + \Gamma_{DD'C'}^{\ \ \ \ \ \  E'}\e_{C E'}),
$$
together with the symmetry requirement
$\Gamma_{DD'CE}=\Gamma_{DD'EC}$,
$\Gamma_{DD'C'E'}=\Gamma_{DD'E'C'}$. These conventions result in the
following expressions for differentiation of spinor \emph{components},
where $\iota^A$ is a two-component spinor field over the manifold:
\begin{eqnarray*}
\nabla_{BB'} \iota^A = \e_{BB'}(\iota^A) + \Gamma_{BB'C}^{\ \ \ \ \ \ A} \iota^{C}, \\
\nabla_{BB'} \iota_A = \e_{BB'}(\iota_A) - \Gamma_{BB'A}^{ \ \ \ \ \ \ C} \iota_{C},
\end{eqnarray*}
and similarly for a primed spinor field. These are the concrete
expressions for the covariant differentiation of spinors using the
connections on $S$ and $S'$ inherited from the Levi-Civita connection,
mentioned in Section \ref{spinorsinneutralsig}. One can extend the above
expressions to multi-component objects in the obvious way, allowing
covariant differentiation of tensors, which agrees with
covariant differentiation using the Levi-Civita connection.

The Riemann tensor has the following spinor decomposition
(\cite{penroserindler}, pg. 236):
\begin{eqnarray} \label{riemann}
R_{abcd} &=& C_{ABCD} \ep_{A'B'} \ep_{C'D'} + \tilde{C}_{A'B'C'D'}
\ep_{AB} \ep_{CD} \nonumber\\ 
&&+ \Phi_{ABC'D'} \ep_{A'B'}\ep_{CD} +
\Phi_{A'B'CD} \ep_{AB}\ep_{C'D'} \nonumber\\&& + 2 \Lambda (\ep_{AC}\ep_{BC}
\ep_{A'C'} \ep_{B'D'} - \ep_{AD} \ep_{BC} \ep_{A'D'} \ep_{B'C'}).
\end{eqnarray}
The Weyl spinors $C_{ABCD}, \tilde{C}_{A'B'C'D'}$ are completely symmetric, and
the traceless Ricci tensor 
$\Phi_{ABC'D'}$ is symmetric on each pair of indices. The
$C,\tilde{C}$ spinors make up the self-dual and anti-self dual parts
of the Weyl tensor. In the language of representation theory, this is
the decomposition of $R_{abcd}$ into irrreducible representations
under the action of $SL(2,\mathbb{R}) \times SL(2,\mathbb{R})$ (with
$\mathbb{R}$ replaced by $\mathbb{C}$ for the holomorphic case).

Note that in $++--$, spinor components are real. For analytic metrics, we can analytically continue which amounts to allowing the spinors to be complex. The remaining calculations in this appendix are valid in both cases.
\vskip5pt
{\bf Integrability of $\alpha$ and $\beta$ surfaces.}
We now show that (\ref{GSF1}) and (\ref{GSF2}) are equivalent to the
fact that the two-plane distributions defined by $o^{A'}$ and
$\iota^A$ are integrable. The leaves are called \as s and \bs s
respectively. The argument is well-known in twistor theory. We will do
the calculation for the $o^{A'}$ case; the $\iota^A$ case is identitical. 

Let $X=\alpha^{A} o^{A'} \e_{AA'}$, $Y=\beta^{A} o^{A'} \e_{AA'}$ be
vector fields, which by definition are in the $\alpha$-planes
determined by $o^{A'}$. Then if they commute we have:
$$
[X,Y]_{AA'} = (f \alpha_A + g \beta_A )o_{A'},
$$
for some functions $f,g$. Multiplying by $o^{A'}$ gives
\begin{equation*}
o^{A'} [X,Y]_{AA'} = o^{A'} (X^{BB'} \nabla_{BB'} Y_{AA'} - Y^{BB'}
\nabla_{BB'} X_{AA'}) = 0.
\end{equation*}
Substituting the spinor expressions for $X^{AA'}$ and $Y^{AA'}$
results in 
$$
o^{A'} o^{B'} \nabla_{BB'} o_{A'} = 0,
$$
which is (\ref{GSF2}), and it is easy to show this is sufficient as
well as necessary.
\vskip5pt
{\bf Twistor distribution and ASD.} 
Locally, the primed spin bundle $S'$ is isomorphic to $M \times
\C^2$. We choose the coordinates on the $\C^2$ to be $\pi^{A'}$ for
$A'=0,1$. This vector bundle has a connection inherited from the
Levi-Civita connection of the metric, and therefore we can find the
horizonal lifts $\et_{AA'}$ of the $\e_{AA'}$, defined by covariantly
constant sections. These lifts are as follows:
$$
\et_{AA'} = \e_{AA'} - \Gamma_{AA'B'}^{\ \ \ \ \ \ C'} \pi^{B'}\frac{\p}{\p \pi^{C'}}.
$$

Using the following formula  (\cite{penroserindler}, pg. 247)
relating curvature quantities to the derivatives of $\Gamma_{AA'C'}^{\
  \ \ \ \ \ D'}$: and the spinor decomposition of the curvature 
(\ref{riemann}) we find 
\begin{eqnarray*}
[\pi^{A'}\et_{AA'},\pi^{B'}\et_{BB'}] &=& (\Gamma_{AA'B}^{\ \ \ \ \ \
  D}-\Gamma_{BA'A}^{\ \ \ \ \ \ D})\pi^{A'} \pi^{B'} \et_{DB'}  \\ && +
\pi^{A'} \pi^{B'} \ep_{AB}\ep^{F'Q'} \tilde{C}_{A'B'E'Q'} \pi^{E'}
 \frac{\p}{\p \pi^{F'}}.
\end{eqnarray*}
One can see from this that if $\tilde{C}_{A'B'C'D'}=0$ then $\pi^{A'}
\et_{AA'}$, $A=0,1$, forms an integrable distribution. The projection
of a leaf of this distribution to $M$ gives an \as. We have
demonstrated that if the metric is anti-self-dual, then given any
point $p \in M$ and an $\alpha$-plane at $p$, there is a unique \as \
through $p$ tangent to this $\alpha$-plane . This was first shown by Penrose
\cite{penrose}, although without using the primed spin bundle. For our
purposes the above formulation will be most useful.

{\bf Proof of Proposition \ref{liftofK}.} 
\label{propproof}
We have the following identity:
\begin{equation} \label{killingidentity}
K^a R_{abcd} = \nabla_b \nabla_c K_d - \frac{1}{2}(\eta_{,b} g_{cd} -
\eta_{,c}g_{bd} + \eta_{,d} g_{bc} ),
\end{equation}
where $\eta$ is the conformal factor appearing in $(\ref{killing})$.
Using the curvature decomposition (\ref{riemann}) to convert this into
spinor form, one can calculate
\begin{eqnarray} 
\label{secondcom}
&&[K^{AA'} \et_{AA'}, \pi^{B'} \et_{BB'}] = (K^{AA'} \Gamma_{AA'B}^{\ \
  \ \ \ \ D} - \psi_B^{\ \ D})L_D\\  
&&- \pi^{B'} (\phi_{B'}^{\ \ A'}
  \ep_B^{\ \ A} + \frac{1}{2} \eta \ep_{B'}^{\ \ A'} \ep_{B}^{\ \
  A})\ \et_{AA'} \nonumber\\
&&+ \pi^{B'} \pi^{E'} ( \e_{BB'} \phi_{E'}^{\ \ F'} -
  \Gamma_{BB'E'}^{\ \ \ \ \ \ \ G'}\phi_{G'}^{\ \ F'} 
+
\Gamma_{BB'G'}^{\ \ \ \ \ \ \ F'} \phi_{E'}^{\ \ G'} - \frac{1}{4}(\e_{BE'}
\eta) \ep^{F'}_{\ \ B'}) \frac{\p}{\p \pi^{F'}}\nonumber.
\end{eqnarray}
We wish to add a vertical term to $K^{AA'} \et_{AA'}$ which will
cancel all the non-$L_A$ terms on the RHS of (\ref{secondcom}). We
don't mind about multiples of the Euler vector field since this gets
quotiented out on projectivizing. A simple calculation shows that
$\tK$ as defined in $(\ref{lielift})$ does the trick.
$\square$

\end{document}